\documentclass{itor_preprint}

\usepackage{natbib}%
\usepackage[figuresright]{rotating}

\usepackage{url}
\usepackage{algorithmic}
\usepackage{algorithm}

\usepackage{amsmath,amsthm}

\usepackage{hyperref} 
\usepackage{cleveref} 
\usepackage{amsfonts} 
\usepackage{bm}       
\usepackage{enumitem} 
\usepackage{multirow} 
\usepackage{comment}  

\usepackage{xcolor}

\begin{document}

\title{\vspace{-1cm}A stochastic integer programming approach to reserve staff scheduling with preferences}

\author[Carl Perreault-Lafleur]{Carl Perreault-Lafleur\affmark{a,$\ast$},  Margarida Carvalho\affmark{a} and Guy Desaulniers\affmark{b}}

\affil{\affmark{a}CIRRELT and Département d’informatique et de recherche opérationnelle, Université de Montréal, Montréal, H3T 1J4, Canada}
\affil{\affmark{b}GERAD and Départment de mathématiques et de génie industriel, Polytechnique Montréal, Montréal, H3T 1J4, Canada}
\email{carl.perreault-lafleur@umontreal.ca [Carl]; carvalho@iro.umontreal.ca [Margarida];\\ guy.desaulniers@gerad.ca [Guy]\vspace{1cm}}

\thanks{\affmark{$\ast$}Author to whom all correspondence should be addressed (e-mail: carl.perreault-lafleur@umontreal.ca).}


\begin{abstract}
    Nowadays, reaching a high level of employee satisfaction in efficient schedules is an important and difficult task faced by companies. We tackle a new variant of the personnel scheduling problem under unknown demand by considering employee satisfaction via endogenous uncertainty depending on the combination of their preferred and received schedules. We address this problem in the context of reserve staff scheduling, an unstudied operational problem from the transit industry. To handle the challenges brought by the two uncertainty sources, regular employee and reserve employee absences, we formulate this problem as a two-stage stochastic integer program with mixed-integer recourse. The first-stage decisions consist in finding the days off of the reserve employees. After the unknown regular employee absences are revealed, the second-stage decisions are to schedule the reserve staff duties. We incorporate reserve employees’ days-off preferences into the model to examine how employee satisfaction may affect their own absence rates.
\end{abstract}


\maketitle

\section{Introduction}\label{sec:Intro}
    Employee absences can have a critical impact on the quality of service offered by companies. This is particularly the case in the transportation sector where the absenteeism is very high due to atypical workdays, which can easily lead to service cancellations. Therefore, some enterprises rely on a specific pool of workers, called extra-boards (XBs) in the transit industry, to mainly cover the regular employees’ absences. In contrast to the known-in-advance absences such as employee vacations, long term absences and open work, the unknown absences that are declared close to the operation day result in uncertainty about when to schedule XBs. It is important that each day, the number of XBs scheduled matches the demand in terms of the shortage of regular employees to avoid service cancellation or XBs idle time. In the remainder of the study, when we mention demand, we are referring to demand in terms of XBs to cover the regular employee absences.
    
    As described in \cite{DESAULNIERS200769}, the planning process of transit industry is commonly divided into three parts: strategical, tactical and operational. The strategical part aims at taking long-term decisions having a direct impact on the quality of service such as the network and route transit design. The tactical phase is reviewed on a seasonal basis, and also concerns the quality of service such as the bus frequencies and timetabling, recruitment targets and employee dispatch by division. Finally, the operational phase seeks to offer the previously proposed service at a minimal cost. To do so, many different problems are solved, the largest two problems being vehicle scheduling and personnel scheduling. This is where the XB (personnel) scheduling problem arises. The XB scheduling process varies around the world. In Europe, the XBs are assigned each possible days-off pattern and duty (the shift of the employees, consisting of the start time, end time, and a potential pause) over time, as part of rotating schedules. In contrast, in North America, the days-off patterns are assigned from scratch to the XBs at the start of each scheduling horizon, and the duties are assigned at the latest the eve of the operation day. In practice, the dispatcher decides on the number of days off to be offered each day of the scheduling horizon. Then, according to their seniority, the XBs take turn choosing which days off they prefer among the ones still available.
    
    Due to ongoing and anticipated labor market shortages, it has become essential for companies to take into account employee preferences to differentiate themselves, recruit new employees and address employee retention. In our study, we follow the North American XB scheduling process and directly assign the days-off pattern to the XBs to better incorporate their preferences in the optimization process. This approach generalizes well to multiple applications that consider employee preferences, such as in nurse scheduling (see \cite{HOLMES76} and \cite{Goodman2009}).

\subsection{Contributions}\label{sec:Cont}

    Our contributions summarize as follows. 
    
    First, we formulate a two-stage stochastic integer program with mixed-integer recourse to model the XB days-off scheduling problem over a finite time horizon. In most cases, this scheduling is done manually by the dispatchers, at best making use of average absence ratios. In our model formulation, we leverage the decomposable aspect of absence scenarios per time period to achieve representational power for uncertainty. This contrasts with the classical method of defining scenarios in stochastic programming, which would require exponentially more scenarios. 
    
    Second, we consider XB preferences in the process to model uncertainty. Social science researches aiming at reducing employee absenteeism have unanimously identified job satisfaction as a key factor influencing the employee's motivation to attend work (\cite{Kehinde11}, \cite{Tasie18}). Personal reasons, at the individual level, have also been identified as an absenteeism cause. For example, some employees might have some strict constraints (due to family responsibilities, health issues, etc.), making them unavailable to work on certain days. These motivate our initial assumption that an XB assigned to preferred days off is less likely to be absent. This assumption not only gives the flexibility to increase the satisfaction level of employees receiving their preferred days off, but also allows to manage the strict schedule constraints some employees might have. To the best of our knowledge, this is the first time that reserve employee preferences are considered in an absence staffing problem.
    
    Third, we expose a new type of problem mixing both exogenous and endogenous uncertainties. In one hand, exogenous uncertainty, i.e., uncertainty independent from the decisions, is the most standard form of uncertainty in stochastic programming and is represented in our problem as the demand. On the other hand, endogenous uncertainty has been little studied and is much more difficult to treat as the stochastic processes are affected by the decisions. Indeed, in our problem, when an XB is assigned to a days-off pattern, he/she can be satisfied, unsatisfied or neutral, which can alter his/her own absence probability on certain of his/her working days. Hence, the uncertainty regarding XB absences depends on the first-stage decisions related to the days-off pattern assignment. In our formulation, we mimic this endogenous uncertainty in a linear fashion.
    
    Our approach is validated with an empirical study based on data from the city of Los Angeles. We show that taking preferences into account leads to significant improvements when compared to not considering preferences: on average, employee social welfare is improved by 37.12\% while also substantially minimizing the cancelled service by 21.14\%. The stochastic programming formulation is partly responsible for theses gains as the social welfare and the cancelled service are improved on average by 7.65\% and 16.63\%, respectively, when comparing with a deterministic formulation. 
    
\subsection{Paper organization}\label{sec:PaperOrg}
    This paper is organised the following way. \Cref{sec:RelLit} surveys the literature relevant to our problem. In \Cref{sec:ProbDesc}, we describe the XB days-off scheduling problem and the proposed mathematical program to model it in \Cref{sec:ModelForm}. We detail the setup of our empirical study in \Cref{sec:EmpStudySetup}, and report the results in \Cref{sec:CompExp}. Finally, \Cref{sec:Conc} draws conclusions and discusses possible directions for future research.
    
\section{Related literature}\label{sec:RelLit}

    In this section, we provide a detailed landscape of the existent literature related to personnel scheduling, uncertainty in scheduling and preference optimization.
    
    \paragraph{\textbf{Personnel scheduling}}
        Historically, personnel scheduling problems have been the subject of numerous studies in many industries: transportation, health care, retail stores, etc. It consists in finding the days off and working days of a set of employees, and the shifts or duties for the working employees. In general and at a very high level, the set of hard constraints is divided into two parts. The first one requires a covering of the demand by the employees. The second one can be seen as a set of work rules to comply to, ensuring the feasibility of the employee schedules. Then, a set of soft constraints comes on top to define which solutions are preferable among the set of valid schedules. Employee preferences, for example, belong to this class of soft constraints. These groups of constraints are easily identified in most of the personnel scheduling problems, e.g., the nurse scheduling problem (\cite{Jafari2015}), the airline and transit crew rostering problems (\cite{Kohl2004} and \cite{Xie2012}), and the retail store workforce scheduling problem (\cite{Chapados2011}). Of course, many variations emerge from one domain to another. For instance, cyclic schedules are often required in the European bus transit industry, employee competency-based schedules are mostly used in the health care industry, and preference-based schedules are increasing in popularity across all domains.
        
        Within the personnel scheduling literature, it is particularly relevant the works on days-off scheduling. To the best of our knowledge, no literature exists on our particular problem of assigning days off to reserve employees. However, general days-off scheduling has been studied. Days-off scheduling consists in finding the right daily number of employees needed to satisfy the daily demand, while ensuring the days-off rules hold. Multiple policies exist for the days off: two days off per week, two consecutive days off per week, four days off every two weeks, etc. This problem alone is not very complex as the number of days-off patterns is usually low and the demand is pre-determined. Its main challenge resides in the modeling of the days-off rules. Variants include consideration of multiple types of employees such as part time and full time (\cite{Emmons97}), employee skills and qualifications for tasks (\cite{ULUSAMSECKINER2007694}), and cyclic schedules (\cite{Emmons91}). After having decided who is working on which days, the shift scheduling occurs, which consists in assigning duties to working employees to satisfy the demand over the course of the day. Shift scheduling alone can have varying complexity depending on the number of possible duties. Again, the employee demand is pre-determined. Embedding days-off with shift scheduling yields the so-called tour scheduling. When the duty start and end times are mostly invariant (such as in manufacturing companies where the duties are often 9AM-5PM), the tour scheduling problem can be tackled directly. Such an example is shown in \cite{BAILEY1985395} which considers only five different start times. However, usually, the days-off and shift scheduling problems must be solved in turn, due to the explosion of possibilities when considering the combinations of days-off patterns and duties. \cite{vanVeldhoven2016} show that although this 2-step decomposition reduces the solution time by 80\% to 90\%, the quality of the solutions is often deteriorated compared to when solving directly the tour scheduling problem. The problem tackled by us is a (stochastic) tour scheduling problem.

    \paragraph{\textbf{Scheduling under uncertainty}}
        \cite{VANDENBERGH2013367} offer a broad review of the solution methods for general personnel scheduling problems. In their review, they point out three uncertainty sources that could be faced. Those are uncertainties related to the demand (what is the workload to accomplish?), the arrival (when does the workload occur?), and the capacity (how many employees can be used?). In the personnel scheduling literature, it is very common to account for pre-determined demand and number of employees. Those uncertainties about the workload and the workforce are thus, most of the time, completely ignored. In these deterministic approaches, the workers are assumed to be always present, whereas absences are unavoidable in practice. The uncertainties related to the workload are also disregarded; when not fixed and known, the workload is estimated via historical data, using forecast and prediction techniques such as machine learning. Although these estimations can be accurate, the stochastic aspect of the demand is left out, leaving the developed models sensitive to when the estimations differ from the real demand. Some approaches exist to account for the disruptions while keeping a deterministic modeling. In this context, the goal is to create schedules that necessitate only few adjustments when disruptions happen to the estimated demand. Examples of works in this line are \cite{Ingels2017,Paias2021}.

        Stochastic approaches, however, are more appropriate due to their ability to incorporate and handle such forms of uncertainty compared to the deterministic approaches. Overall, although most of the personnel scheduling problems are approached in a deterministic way, some papers follow the stochastic path. This is the case of \cite{Kim15}, who formulate an integrated nurse staffing and scheduling problem with unknown demand as a two-stage stochastic integer program. They schedule the nurses for 3-month long horizons, while respecting rules on the nurse-to-patient ratios and not knowing in advance the number of patients in the hospital. Their recourse action is to add or cancel shifts, with associated costs for each. Penalty costs are also inferred for overstaffing and understaffing. They demonstrate computationally the efficiency of their model on the cost savings, and that these savings increase with the precision of the demand forecasts. In our model, the uncertainty, i.e., the stochastic process, will be in personnel absences.
        
        In stochastic optimization, most of the uncertainty we encounter is modeled as exogenous uncertainty. That is, the decision variables do not affect the uncertainty distribution. However, with endogenous uncertainty, the decisions can influence the probability distribution. \cite{Goel2006} and \cite{Li21} distinguish two types of endogenous uncertainty. In type I endogenous uncertainty, decisions influence the parameter realizations by altering the underlying probability distributions for the uncertain parameters. In type II, the decisions influence the parameter realizations by affecting the time at which we observe these realizations. In our setting, only the type I is of interest since we hypothesize that the XBs' absence rates vary according to their satisfaction with the assigned days-off patterns. However, to the best of our knowledge, there is no work on scheduling considering type I endogenous uncertainty.

    \paragraph{\textbf{Individual preference optimization}}
        Individual preferences are now fundamental in the design of schedules. More and more personnel scheduling papers incorporate them in their models. In the literature, we can find two different approaches for considering preferences. The first one, the most common, is to account for the preferences directly in the model and maximize for overall satisfaction. This method is used in \cite{BADRI1998303}, \cite{Jafari2015} and \cite{BARD2005510}, to name just a few. The second one consists in first creating schedules without preferences, and then performing an auction where the employees bid on their preferred schedules as in \cite{Grano2008}. In general, few attention is given to the fairness in the attribution of the schedules to employees. Multiple optimal solutions that assign different schedules to different employees might exist, with some employees being more satisfied than others about their schedules. This demonstrates that even if the employee satisfaction about the schedules is globally high, unfairness can arise from an egalitarian perspective. To counter this effect, \cite{BADRI1998303} design the penalty of not respecting a preference as an increasing function of the number and severity (according to some rules) of preference violations. As explained in \Cref{sec:Cont}, one of the reasons to consider employee preferences in schedules is to increase job satisfaction, which can in turn decrease absenteeism. The inverse is also, and perhaps even more, true: \cite{deBoer2002} state that job dissatisfaction, together with stress, are the two explanations for absenteeism. In the mathematical programs designed to solve preference-aware personnel scheduling problems, keeping track of employee satisfaction at the individual level based on the realization of their preferences on the assigned schedules becomes challenging because it depends on decision variables.

\section{Problem description}\label{sec:ProbDesc}
\begin{figure}[h!]
        \centering
        \includegraphics[width=0.66\textwidth]{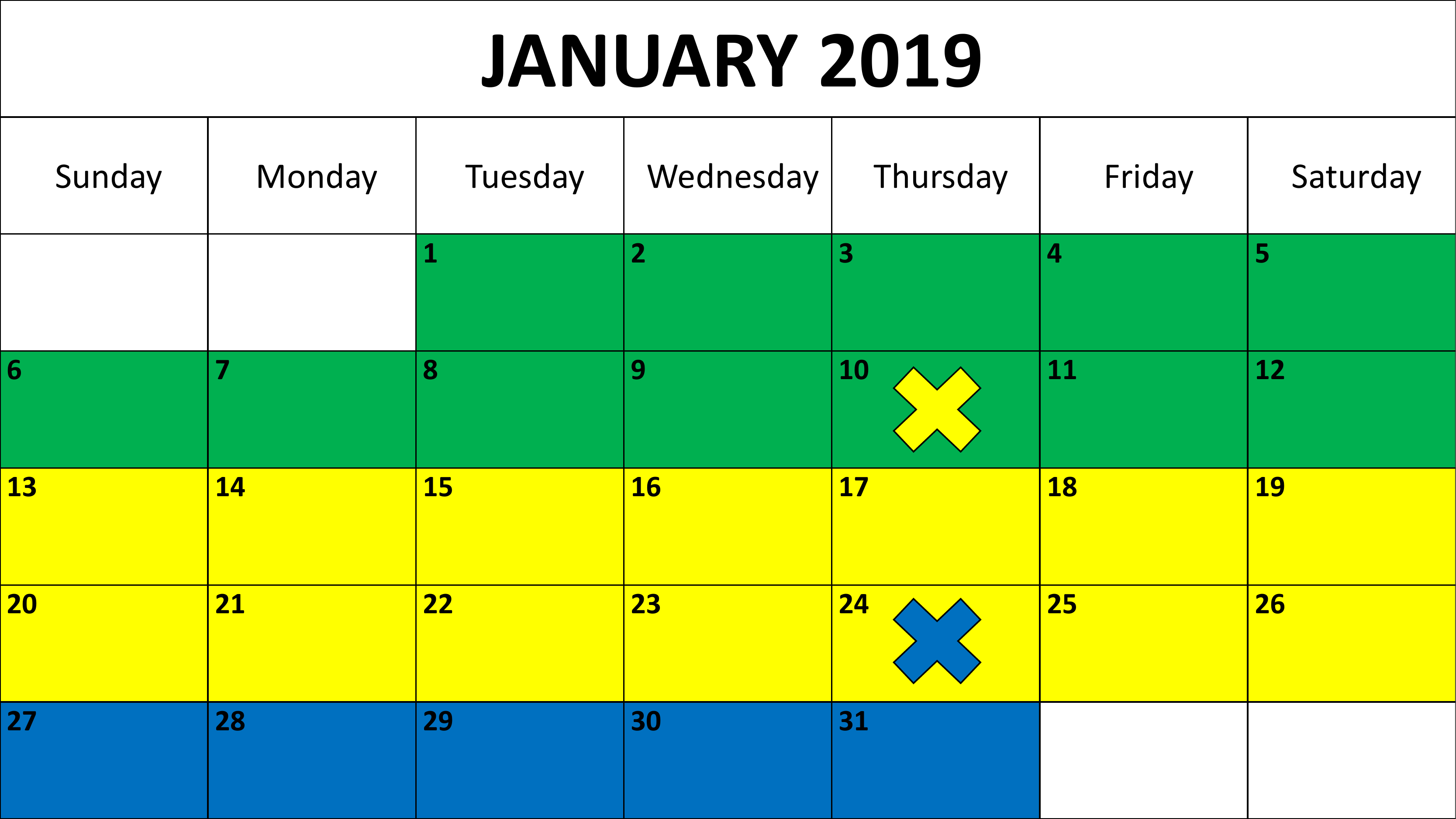}
        \caption{The January 2019 XB days-off assignment calendar in Los Angeles. The crosses represent the days where the assignment of days-off patterns to XBs must be decided, and the same-color days represent the corresponding planning horizon.}
        \label{fig:calendar}
    \end{figure}

    The solving frequency and horizon length of the XB days-off scheduling problem vary from one transit company to another and range from 2 to 4 weeks. In Los Angeles, the problem is solved bi-weekly in all the divisions, for a planning horizon $J$ of 2 weeks long ($|J|=14$ days). We can see an example of this process in \Cref{fig:calendar}. Next, we start by detailing the first-stage decisions of our planning problem, followed by a description of the XBs demand including the associated uncertainty and, finally, the recourse decisions.
    
    The goal of the XB days-off scheduling problem is to obtain a work schedule for each XB in the set $E$ of XBs, first consisting of his/her days off and working days for each of the $|J|$ days of the planning horizon. Then, in a second time, given a demand scenario, the set of duties that the planned-to-work XBs will receive must be found, but these duties do not need to be individually assigned. These two types of decisions correspond respectively to the first stage and the second stage (recourse) of the stochastic problem to be formulated. Throughout the rest of the paper, we will refer to those work schedules, i.e., first-stage decisions, as days-off patterns. We define $P$ to be the set of days-off patterns. In Los Angeles, the XBs must be assigned to two consecutive days off per week that must be the same for each week of the horizon, making a total of 7 possible days-off patterns: Sunday-Monday, Monday-Tuesday,..., Saturday-Sunday.
    
    Two types of demand coexist in the problem: known and unknown at the moment of assigning the days-off patterns. The unknown demand corresponds to a stochastic variable and we represent its distribution using scenarios. Let us clarify in detail the demands for XBs:
    \begin{itemize}
        \item For each day $j\in J$, $o_j$ represents the number of unassigned duties from open work and known-in-advance absences, vacations and long-term absences. One XB can fulfill exactly one such duty, on a given day. Covering this known demand is formulated with hard constraints to ensure a reasonable number of XBs is available to cover absences.
        \item For each day $j\in J$, $S_j$ represents the set of daily scenarios. The uncertainty from each scenario $s\in S_j$ is made available through the parameters $d^s_{t,j}$ which depicts the number of unknown absences (absences that were declared after the assignment of the days-off patterns) at time period $t\in T$ on day $j$, where $T$ is the set of time periods. For example, $T$ can be chosen to be hourly time periods, in which case $|T|=24$.
    \end{itemize}
    In practice, after the days off have been assigned, the XBs know on which days they are scheduled to work. XBs working on known-in-advance absences and open work are also assigned their duties. However, for the XBs working on unknown absences, the duties that they will inherit are still to be defined. At the latest, those duties are assigned to these XBs on the eve of each working day. At that moment, most of the unknown absences are declared and the demand in terms of XB working hours to cover the regular employees' absences is almost fully known. This motivates the modeling of the problem as a two-stage stochastic program, where the recourse decisions correspond to determining the XBs' work shifts (duties) for each absence scenario. Concretely, given a scenario $s\in S_j$ for day $j \in J$, each XB working on fulfilling unknown absences should be assigned a duty $w$ from the set of duties $W$. The set of duties depend on the labor union agreement. Each duty has a specific cost $c^w_2$, depending on the work, pause and overtime hours it contains.
    
    In addition to the regular employees' absences, the XBs' absences have to be considered in the problem. Indeed, some XBs can also be absent during their own planned duties, creating a shortage in the number of XBs scheduled to cover the regular employees' absences. The probability of absence of XB $e\in E$ given that he/she is assigned to work on day $j\in J$ is denoted by $q_{e,j}$. Our assumption is that this probability is linked to the XB's happiness about his/her assigned days-off pattern. If the XB is off on his/her most preferred days, we assume that his/her absence probability will be low during the assigned working days. Inversely, if the XB is planned to work on his/her most preferred days, the absence probability will be higher for these days. More details about this process are given in \Cref{sec:ParamSettings}. Each XB $e\in E$ expresses his/her preferences about the $7$ days-off patterns in $P$ by ranking them from 1 to 7 (7 being the most preferred, 1 being the least preferred) into the preference score parameters $s_{p,e}$ for $p\in P$.
    
    Once the duties have been assigned for each scenario, it becomes possible to compute whether understaffing (which results in cancelled service) or overstaffing occurs at each time period $t\in T$, for each day $j\in J$. 
    The objective function corresponds to the expected costs for cancelled service and duty costs throughout the entire planning horizon. The expectation is taken with respect to the probability distribution of each scenario. The probability that a scenario $s\in S_j$ occurs is denoted by $\alpha^s_j$, and the probabilities of daily scenarios sum to 1. A penalty of $c_1$ is issued for each period of understaffing, and a duty cost $c^w_2$ for duty $w$ must be paid for each XB working on unknown absences. Duty costs of XBs working on known-in-advance absences and open work are omitted because they are constant (recall that we impose this demand to be covered). Additionally, the objective function also includes the negative sum of the XB preferences, associated with a social welfare benefit $c_3$. Thus, our goal is to minimize the described objective function. It is important to note that minimizing the cancelled service and maximizing the XB preferences (social welfare) is closely related and can be jointly optimized. Indeed, increasing the XBs' satisfaction by assigning days-off patterns they like will make them less absent, creating more flexibility to cover the regular employee absences and avoid service cancellation.
    No cost is assigned to overstaffing because there is already an implicit cost to it: precious work time from the XBs is wasted by not covering unknown absences.
    
\section{Mathematical formulations}\label{sec:ModelForm}
    In this section, we formalize the XB days-off scheduling problem as a two-stage stochastic program with mixed-integer recourse. We consider two cases, namely, with and without preferences (\Cref{sec:ModelFormPref,sec:ModelFormNoPref} respectively). The notation used in these models is listed in \Cref{tab:Notation}. 
    
    \begin{table}[h!]
    \centering
    \caption{Notation}
    \label{tab:Notation}
    \begin{tabular*}{\hsize}{@{}@{\extracolsep{\fill}}ccl@{}}
    \hline
    \textbf{Type}      & \textbf{Notation} & \multicolumn{1}{c}{\textbf{Description}} \\ \hline
    \multirow{6}{*}{\textbf{Sets}} & $J$   & Set of days in the planning horizon. $|J|$ is a multiple of 7. \\
     & $S_j$ & Set of daily scenarios on day $j$ \\
     & $E$   & Set of XBs \\
     & $P$   & Set of possible days-off patterns \\
     & $T$   & Set of time periods \\
     & $W$   & Set of duties \\ \hline
    \textbf{$\bm{1^\text{st}}$-stage variables} &
      $x_{p,e}$ &
      Takes value 1 if pattern $p$ is assigned to XB $e$ and 0 otherwise \\ \hline
    \multirow{3}{*}{\textbf{$\bm{2^\text{nd}}$-stage variables}} 
     & $y^s_{t,j}$ & \begin{tabular}[c]{@{}l@{}}Number of understaffing hours (cancelled service) at time $t$ of  day $j$,\\ under scenario $s$\end{tabular} \\
     & $z^s_{t,j}$ & Number of overstaffing hours at time $t$ of day $j$, under scenario $s$ \\
     & $v^s_{w,j}$ & Number of XB assigned to duty $w$ on day $j$, under scenario $s$ \\ \hline
    \multirow{7}{*}{\textbf{\begin{tabular}[c]{@{}c@{}}Parameters in the \\ constraints\end{tabular}}} &
      $q_{e,j}$ & Absence probability of XB $e$ if working on day $j$ \\
     & $r_{p,j}$ & Takes value 1 if pattern $p$ plans work on day $j$ and 0 otherwise \\
     & $b_{p,j}$ & Takes value 1 if pattern $p$ plans off on day $j$ and 0 otherwise \\
     & $o_j$ & \begin{tabular}[c]{@{}l@{}}Number of duties to be fulfilled on day $j$ to cover known-in-advance \\ absences and open work\end{tabular} \\
     & $a_{w,t}$ & Takes value 1 if duty $w$ plans work at time $t$ and 0 otherwise \\
     & $d^s_{t,j}$ & Number of unknown absences at time $t$ on day $j$, under scenario $s$ \\
     & $\varepsilon$ & Numerical tolerance parameter. A value of $1\mathrm{e}{-15}$ is used. \\ \hline
    \multirow{5}{*}{\textbf{\begin{tabular}[c]{@{}c@{}}Parameters in the \\ objective function\end{tabular}}} 
     & $\alpha^s_j$ & Probability of scenario $s$ on day $j$. \\
     & $c_1$ & Cancelled service cost per time period \\
     & $c^w_2$ & Cost of duty $w$ \\
     & $c_3$ & Social welfare benefit \\
     & $s_{p,e}$ & Preference score of XB $e$ for pattern $p$. This is a number between 1 and 7. \\ \hline
    \end{tabular*}
    \end{table}

\subsection{Formulation for the case with preferences}\label{sec:ModelFormPref}

    The extensive form of the mixed-integer stochastic program for the case with preferences expresses as follows:
    
    \begin{eqnarray}
         \hspace{-0.85cm} \min \quad \sum_{j\in J} \sum_{s\in S_j} \alpha_j^{s} \left(\sum_{t\in T} c_1 y_{t,j}^{s} + \sum_{w\in W} c_2^w  v_{w,j}^{s} \right) - c_3 \sum_{e\in E}\sum_{p\in P}s_{p,e} x_{p,e} \label{eq:ObjFun} \\
         \hspace{-0.85cm} \mbox{subject to:} \hspace{7.75cm} \sum_{p\in P} x_{p,e} = 1,  && \hspace{-0.4cm} \forall e\in E \label{eq:OnePatternPerXB} \\
         \hspace{-0.85cm} \sum_{e\in E}(1-q_{e,j}) \sum_{p\in P}r_{p,j} x_{p,e} \ge o_j,  && \hspace{-0.4cm} \forall j\in J \label{eq:KnownDemandSatisfied} \\
         \hspace{-0.85cm} \sum_{w\in W} v_{w,j}^{s} \le |E| - o_j - \sum_{e\in E} \sum_{p\in P} b_{p,j} x_{p,e} - \sum_{e\in E}q_{e,j} \sum_{p\in P}r_{p,j} x_{p,e}, && \hspace{-0.4cm} \forall j \in J, s \in S_j \label{eq:UpperBoundNumberduties} \\
         \hspace{-0.85cm} \sum_{w\in W} v_{w,j}^{s} \ge |E| - o_j - \sum_{e\in E} \sum_{p\in P} b_{p,j} x_{p,e} - \sum_{e\in E}q_{e,j} \sum_{p\in P}r_{p,j} x_{p,e} - (1-\varepsilon), && \hspace{-0.4cm} \forall j \in J, s \in S_j \label{eq:LowerBoundNumberduties} \\
         \hspace{-0.85cm} \sum_{w\in W} a_{w,t} v_{w,j}^{s} + y_{t,j}^{s} - z_{t,j}^{s} = d_{t,j}^{s}, && \hspace{-0.4cm} \forall j \in J, s \in S_j, t\in T \label{eq:DemandEquation} \\
         \hspace{-0.85cm} x_{p,e}\in \{ 0, 1\}, && \hspace{-0.4cm} \forall e \in E, p \in P \label{eq:BinaryX} \\
         \hspace{-0.85cm} y_{t,j}^{s}, z_{t,j}^{s}\in \mathbb{N}, && \hspace{-0.4cm} \forall j \in J, s \in S_j, t\in T \label{eq:BinaryYZ} \\
         \hspace{-0.85cm} v_{w,j}^{s} \in \mathbb{N}, && \hspace{-0.4cm} \forall j \in J, s \in S_j, w\in W. \label{eq:BinaryV}
    \end{eqnarray}
    Constraints (\ref{eq:OnePatternPerXB}) ensure that each XB is assigned exactly one days-off pattern. In (\ref{eq:KnownDemandSatisfied}), we require that for each day, the number of scheduled-to-work and present XBs is sufficient to cover the known-in-advance absences and open work. Indeed, $1-q_{e,j}$ represents the presence probability of XB $e$ on day $j$. Summing over all the scheduled-to-work XBs gives the expected number of scheduled-to-work and present XBs. It is important to note that in this model, we account for anonymous XB absences: the number of XB absences is obtained by summing the individual XB absence probabilities, but we do not infer absences of specific XBs. Constraints (\ref{eq:UpperBoundNumberduties})-(\ref{eq:LowerBoundNumberduties}) and (\ref{eq:BinaryV}) can be combined to obtain:
    
    \begin{equation} \label{eq:Sumduties}
     \sum_{w\in W} v_{w,j}^{s} = |E| - o_j - \sum_{e\in E} \sum_{p\in P} b_{p,j} x_{p,e} - \left\lceil \sum_{e\in E}q_{e,j} \sum_{p\in P}r_{p,j} x_{p,e} \right\rceil, \forall j \in J, s \in S_j.
    \end{equation}
    In other words, for each day and demand scenario, the total number of duties we can use to cover unknown absences is the total number of XBs from which we subtract the ones working on known-in-advance absences and open work, the ones that are off, and the ones that were scheduled-to-work but absent. Splitting (\ref{eq:Sumduties}) into two sets of constraints in the mixed-integer program (MIP) is necessary to account for the integrality constraints on the $v^s_{w,j}$ variables. Finally, (\ref{eq:DemandEquation}) ensure that for each day, scenario and time period, the demand equation is fulfilled: the number of XBs working on unknown absences plus the understaffing minus the overstaffing is equal to the unknown absences.
    
    The objective function (\ref{eq:ObjFun}) contains three terms. As described in \Cref{sec:ProbDesc}, the first two, related to cancelled service and duty costs, are scenario-dependent, and so we account for their expected value over the scenario's probability distribution. The last term represents the sum of the XB preferences related to their assigned days-off patterns, often referenced to as the social welfare in economics.
    
    Remark that our mathematical program corresponds to a two-stage stochastic program although not explicitly formulated as such; we have written directly its extensive form by explicitly considering all scenarios. Thus, the recourse problem for the scenario realization $s$ of the associated two-stage stochastic problem corresponds to the optimization with the $x_{p,e}$ fixed and restricted to the variables with superscript $s$.
    
    We now wish to highlight three simplifications within our model:
    \begin{enumerate}
         \item We assume that the uncertainty is entirely revealed at the moment of deciding the XBs' duties. As absences can be declared during the operation day, a multi-stage formulation would be needed to correctly take care of this uncertainty. The latter would greatly complexify the model, thereby our choice for this approximation.
         \item Whereas in practice the XBs inherit from already existing duties (initially assigned to absent drivers, or that were not yet assigned), we do not follow this process for unknown absences. We assume for simplicity that XBs can instantly change from an initial duty to another. In practice, this is not possible as the relieves (places dispatched along duties where a bus has a driver switch and that allows drivers to make a line transfer or take a pause) take a certain time. In our formulation, we aim for situations where the number of XBs working on unknown absences is the same as the total number of unknown absences per time period, even though in practice a larger number of XBs could be necessary to cover the demand, by considering the relieves. On the other hand, by re-arranging and re-optimizing the set of initial duties from unknown absences on the operational day's eve, it would be possible to generate a number of duties that is smaller than the number of initial duties from the absent drivers, and that would still cover the same service. This time, a smaller number of XBs could suffice to cover the demand. Hence, overall, the estimate we make should be close to the real number of XBs necessary to cover the unknown absences.
         \item We ignore two work rules involving duty combinations such as the minimal rest time between two consecutive duties and the maximal weekly work time. This approximation is possible because the duty assignment is only simulated through the second-stage, and no duties are actually assigned at this point.
    \end{enumerate}

    We note that these three approximations are performed only on the recourse problem. After the first-stage decisions have been taken and closer to the operation day, it will be possible to take more precise recourse decisions with more accurate information, and respecting all the work rules.
    
    As mentioned in \Cref{sec:Cont}, we leverage on the separable aspect of our formulation to define daily scenarios instead of horizon-long scenarios like we traditionally have in stochastic programming literature. In order to obtain the same representational power brought by $n\cdot |J|$ daily scenarios ($n=|S_1|=\ldots=|S_{|J|}|$) in our formulation, we need $n^{|J|}$ horizon-long scenarios. This improvement is crucial in our setting as the time to solve the extensive form (also known as deterministic equivalent) of a two-stage stochastic program is typically dependent on the size of the formulation. Thus, a formulation with a large number of scenarios and, consequently, a large number of second-stage variables is expected to demand long computational times.
    
    Lastly, a very important aspect of this model is how the endogenous uncertainty resulting from the satisfaction of XB $e$ about his/her assigned days-off pattern $p$ is modeled linearly on $x_{p,e}$. Since the variables $x_{p,e}$ are binary and constraints (\ref{eq:OnePatternPerXB}) ensure that one pattern is assigned to XB $e$, we are able to represent in a linear way the XBs' absence probabilities conditional on whether the XB is working. 
    
    There are several procedures to generate values for the parameters $q_{e,j}$ and we propose one in \Cref{sec:ParamSettings}. A general principle that should always hold is that the absence probability should be at its highest (lowest) on the days off that are the most (least) desired by an XB.

\subsection{Formulation for the case without preferences}\label{sec:ModelFormNoPref}

    The XB days-off scheduling problem (\ref{eq:ObjFun})-(\ref{eq:BinaryV}) considers by default the XB preferences. The expected number of absent and scheduled-to-work XBs is obtained by summing up the individual absence probabilities $q_{e,j}$ of the scheduled-to-work XBs. We propose in this section some modifications to (\ref{eq:ObjFun})-(\ref{eq:BinaryV}) that would lead to a model not taking into account the XB preferences. Only two changes are required:
    \begin{itemize}
         \item We redefine the parameters $q_{e,j}$ used in (\ref{eq:KnownDemandSatisfied})-(\ref{eq:LowerBoundNumberduties}) as $q_j$, making the absence probabilities equal for all XBs, and independent from the assigned patterns.
         \item We cancel the social welfare benefit in the objective function by forcing $c_3=0$.
    \end{itemize}
    These two changes lead to solutions uninfluenced by the XB preferences. In the empirical study, we will refer to this special case formulation when needed.

\section{Empirical study setup}\label{sec:EmpStudySetup}
    In this section, we explain how we transform the input data to create problem instances, and how we evaluate the quality of solutions.

\subsection{Input data}\label{sec:InputData}
    In Los Angeles, the XB days-off assignment process (also called markup) happens every other Thursday for the planning horizon starting the following Sunday (3 days later) and ending the Saturday two weeks later (see \Cref{fig:calendar}). The data we have access to is from 10 out of the 11 bus divisions of the city and concern years 2018 and 2019 (except specified otherwise). The data is divided into four parts:
    \begin{enumerate}[label=(\alph*)]
         \item \label{item:dataA} For each division and each day of the 2 years, the counts of absences declared before markup (known-in-advance absences) and the total daily absence durations declared after markup (unknown absences). For example, on January $1^\text{st}$ 2018, there were 7 known-in-advance absences, and a total of 88h35min of unknown absences for division 1.
         \item \label{item:dataB} For each division and each day of the 2 years, the distribution of the total daily unknown absence durations per hourly time periods, that we will refer to as absence shape from now on. Following the example above, on the 88h35min of unknown absences, 0h0min take place from 0AM to 4AM, 1h9min take place between 4AM and 5AM, 2h46min between 5AM and 6AM and so on.
         \item \label{item:dataC} A list of XB days-off pattern preferences from one specific markup in 2019.
         \item \label{item:dataD} For each division and each day of 2019, the XB absence rates.
    \end{enumerate}
    Any information specific to the drivers was disregarded from the data. The counts of known-in-advance absences from Data \ref{item:dataA} are directly used for the values of the $o_j$ parameters. The two years of data from Data \ref{item:dataA}-\ref{item:dataB} were used to fit distributions to generate scenarios, and the year 2019 was chosen to create problem instances. During the year 2019, 25 complete two-week long planning horizons starting on Sundays were available, making a total of 250 planning horizons across the 10 divisions.

\subsection{Scenario generation}\label{sec:ScenarioGen}
    In this section, we describe how we generate daily unknown absence scenarios from the data in order to determine values of the parameters $d^s_{t,j}$. The creation of a daily scenario consists of three parts. We start by generating a total daily unknown absence duration (scalar). An example would be a total of 240h of unknown absence on a day. We continue by generating a daily unknown absence shape (vector of length $|T|$ whose values sum to 1). For instance, the shape $(1/24, 1/24, ..., 1/24)$ implies that the unknown absences are distributed evenly over the day. Finally, we combine these two results to obtain the absence scenario. Following our example, we would generate the scenario $(10, 10, ..., 10)$ where there are 10 unknown absence hours during each hour of the day. The process we describe below is for one division, but is repeated for each of the 10 divisions to capture their specific characteristics. 
    
    \begin{enumerate}
         \item To generate the daily unknown absence durations, we consider the set of $2\cdot365=730$ values from the total daily unknown absence durations over the two years from Data \ref{item:dataA}. We find the best-fit distribution $\mathcal{D}(\mu)$ on these values, where $\mathcal{D}$ represents the distribution type (for example a Gaussian) and $\mu$ its parameters (for example mean and variance). We then split this set of values into 7 subsets based on their associated days of the week, and compute the best-fit parameters $\mu_i$ of $\mathcal{D}$ on these subsets $i\in\{1,2,...,7\}$. We obtain, for each division, a globally best-fit distribution with 7 different days-of-the-week-specific best-fit parameters $\mathcal{D}(\mu_1), \mathcal{D}(\mu_2), ..., \mathcal{D}(\mu_7)$ that are independent from the planning horizons. Now, to generate $l$ daily unknown absence durations for a specific day of the week $i$ from any planning horizon, we sample $l$ values from $\mathcal{D}(\mu_i)$.
         \item To generate the daily unknown absence shapes for a specific day of a planning horizon, we simply pick the $k$ daily unknown absence shapes of the same day of the week from the $k$ weeks preceding immediately the horizon from Data \ref{item:dataB}.
         \item To obtain $l\cdot k$ daily scenarios for a specific day $j$ from a planning horizon, we generate $l$ total daily unknown absence durations and $k$ unknown absence shapes as described above. Each of the possible duration and shape pairs is combined to create a scenario that populates $d^s_{t,j}$ for $s\in S_j = \{1,2,...,l\cdot k\}$.
    \end{enumerate}
    
\begin{figure}[h!]
        \centering
        \includegraphics[width=0.48\textwidth]{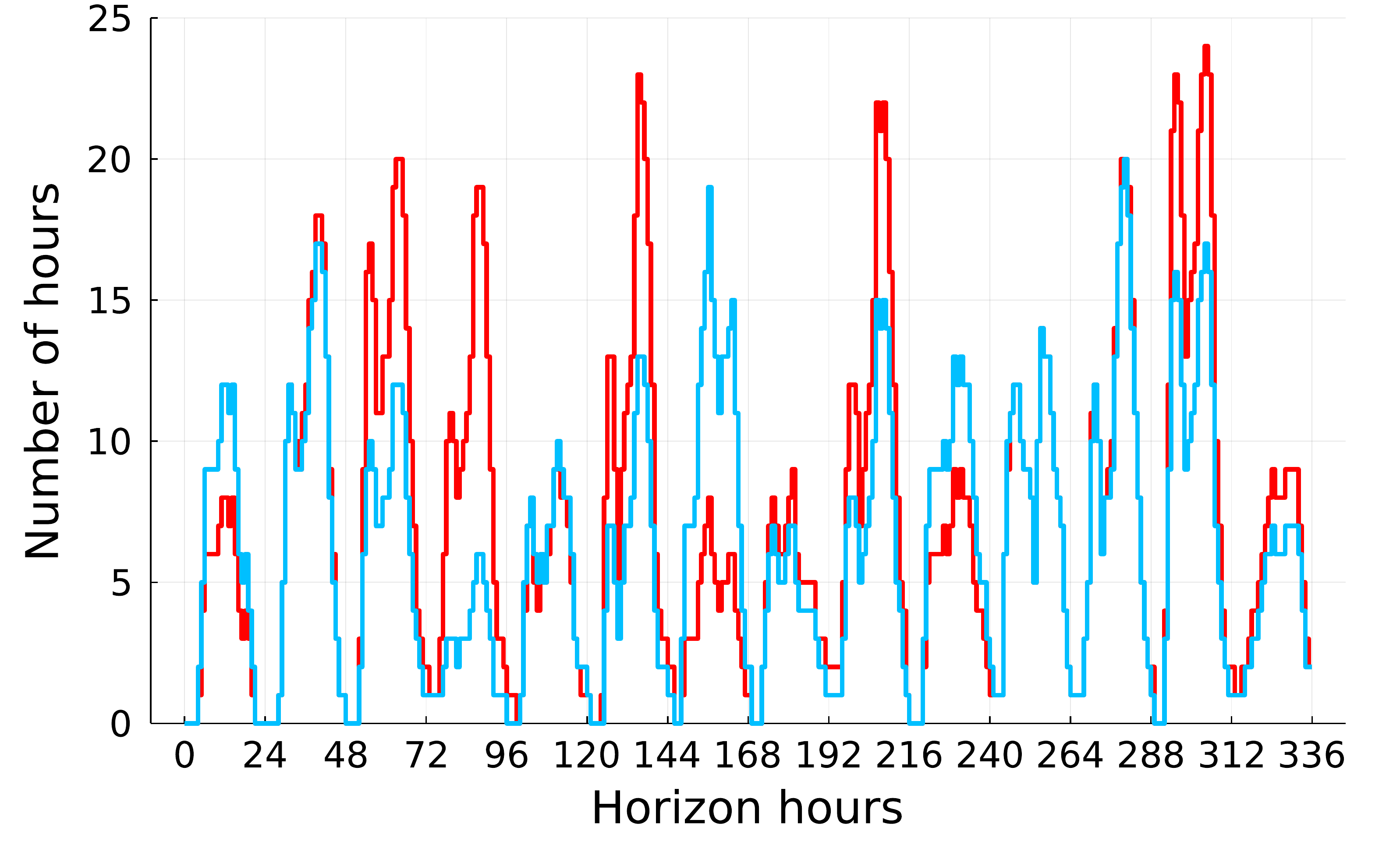}
        \includegraphics[width=0.48\textwidth]{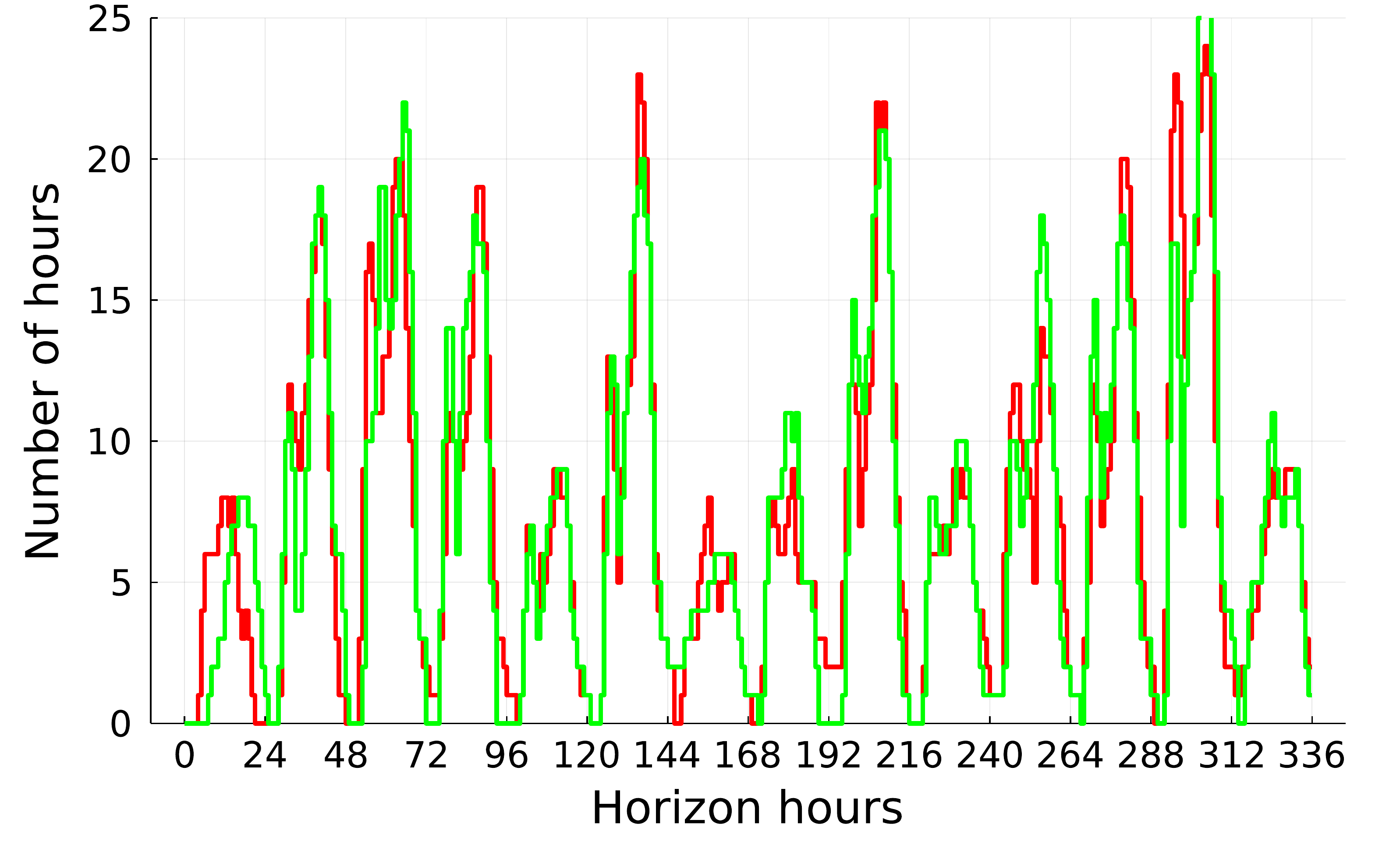}
        \caption{Example of daily scenarios where the blue and red ones use the same absence shapes but different daily absence durations, and the green ones use the same daily absence durations as the red ones but different absence shapes. The daily scenarios are aggregated to form the 14 days horizon.}
        \label{fig:ScenarioGen}
    \end{figure}

    An example of scenarios is provided in \Cref{fig:ScenarioGen}.
    
\subsection{Preference samplings}\label{sec:PrefSamplings}
\begin{figure}[h!]
        \centering
        \includegraphics[width=0.48\textwidth]{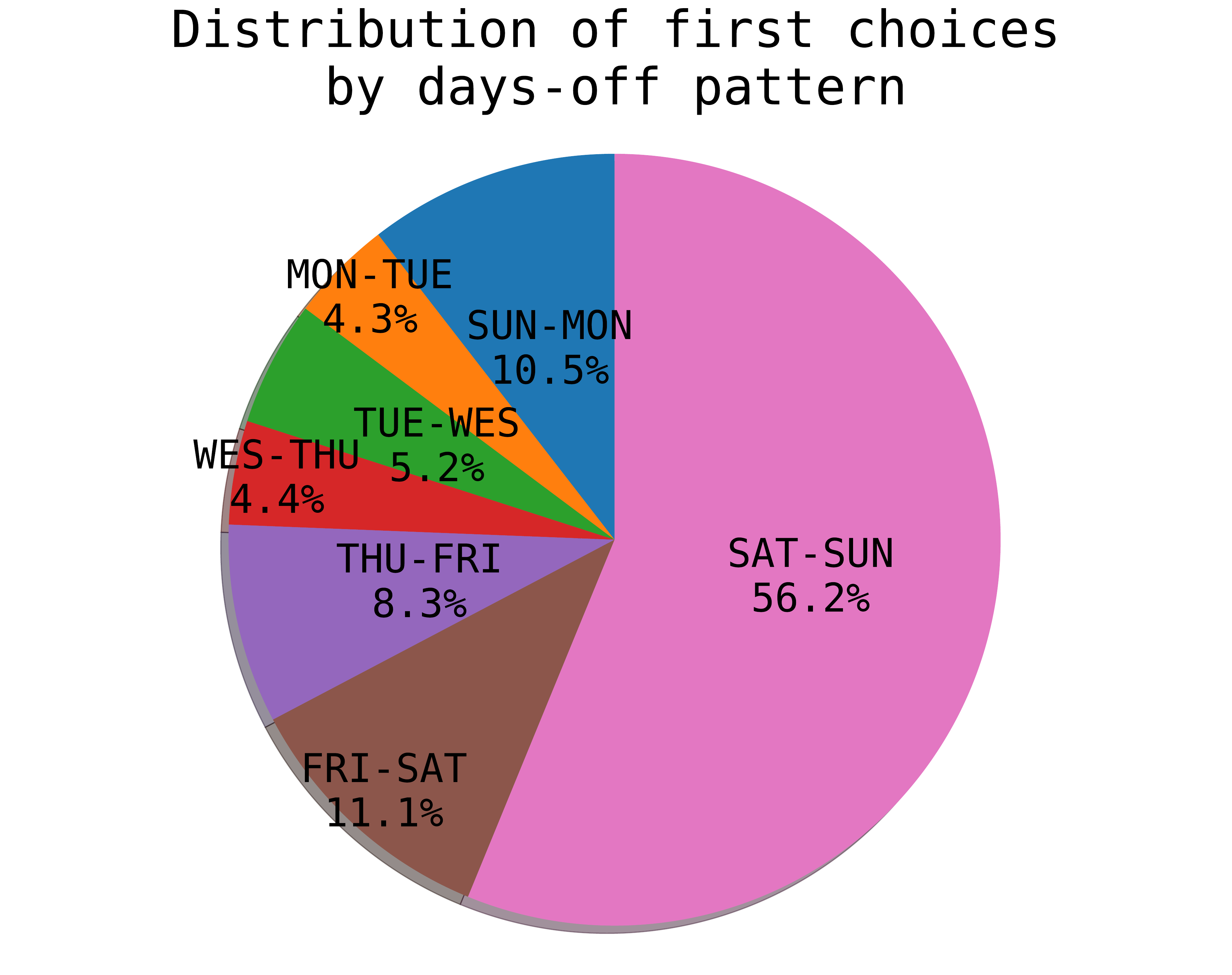}
        \includegraphics[width=0.48\textwidth]{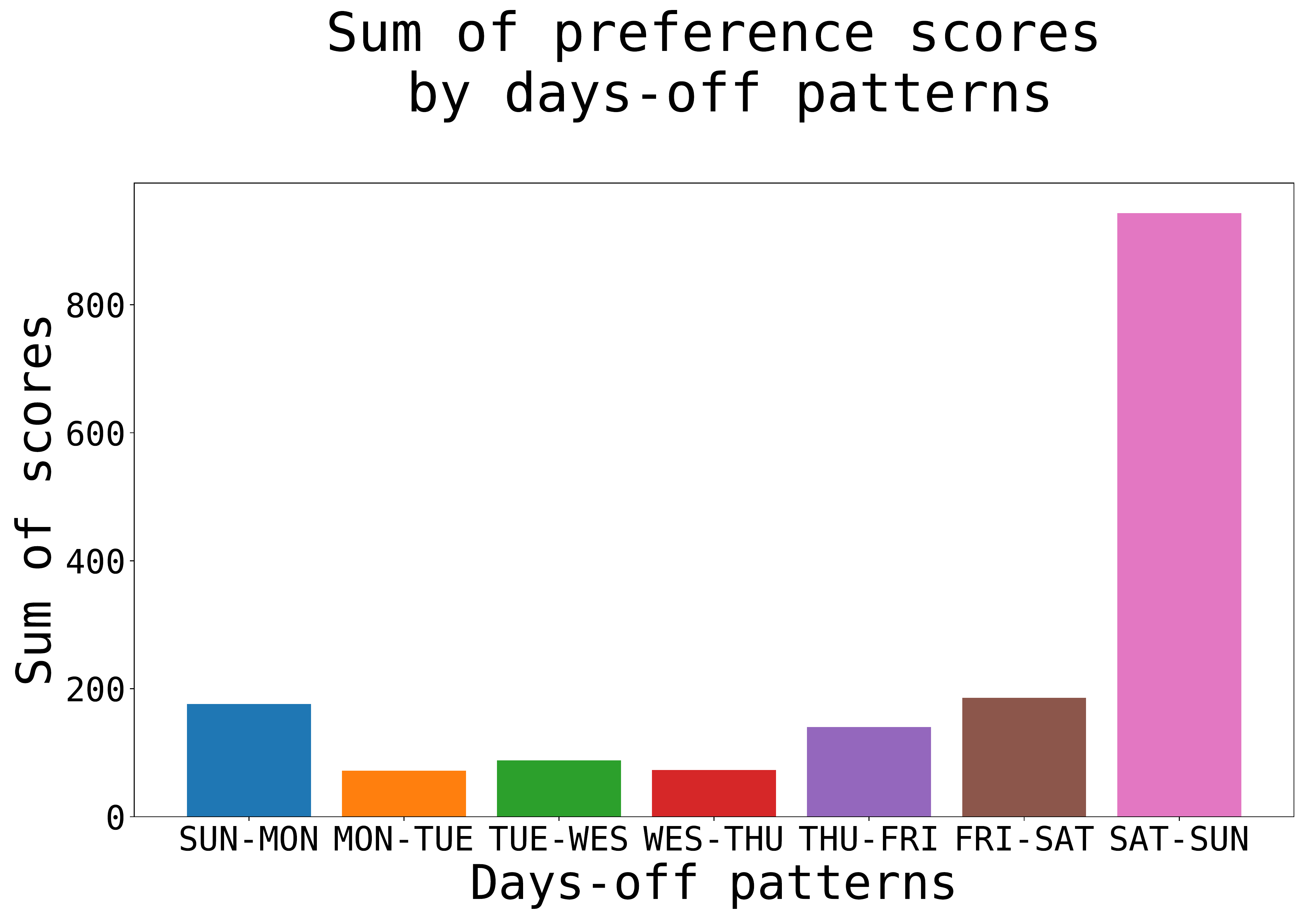}
        \caption{Statistics about the XB days-off patterns preferences.}
        \label{fig:DaysOffPatternsPrefs}
    \end{figure}

    As defined in \Cref{sec:ProbDesc}, the XB days-off preferences are rankings from 7 to 1  (7 being the most preferred, 1 being the least preferred) given to the 7 possible days-off patterns, making a total of $7!=5040$ possible rankings. Data \ref{item:dataC} contains a total of $1678$ days-off rankings from XBs across all divisions. Multiple samplings of days-off preferences from the empirical distribution of these rankings were created to populate the parameters $s_{p,e}$, so that the results we present in \Cref{sec:CompExp} could be averaged over multiple problem instances with different XB days-off preferences. Some statistics about Data \ref{item:dataC} are shown in \Cref{fig:DaysOffPatternsPrefs}.
    \newline

     We can see that the days-off pattern Sat-Sun is by far the most popular, followed by the ones that contain one weekend day (Fri-Sat and Sun-Mon). In fact, the most popular days-off pattern ranking is Sat-Sun (7), Sun-Mon (6), Fri-Sat (5), Thu-Fri (4), Wed-Thu (3), Tue-Wed (2), Mon-Tue (1), and occurs 159 times out of the 1678 rankings. The next 10 most popular rankings also start with Sat-Sun. It is interesting to note that in the transit industry, the number and frequency of buses are much lower on weekends than during the week. Therefore, less drivers are required on these days, which in turn leads to less absences to cover for the XBs. Hence, there is an alignment between preferred days off and the demand for drivers, anticipating an ease in the assignment of days-off patterns that include weekend days. An analogous conclusion can be obtained for the least popular days-off which are those of the week (Mon-Tue, Tue-Wed, Wed-Thu).
    
\subsection{XB absence probabilities}\label{sec:XBAbsProbs}
    The special case formulation where the XB preferences are not considered has a straightforward way of setting the XB absence probability parameters. For each division and day $j\in J$, we simply populate $q_j$ with the average divisional XB absence rate for 2019 of the same day of the week from Data \ref{item:dataD}. Therefore, each day, all the working XBs will have the same absence rate.
    
    When considering the XB preferences, for each division, the historical absence rates from Data \ref{item:dataD} are split by days of the week and the low and high whiskers, $\sigma^\text{low}_i$ and $\sigma^\text{high}_i$, are identified for each day of the week $i\in \{1,2,...,7\}$. The whiskers represent the smallest and largest points from a dataset within 1.5 times the inter-quartile range. Then, for each day of the week $i$, we divide the interval $[\sigma^\text{low}_i, \sigma^\text{high}_i]$ into seven equally distant points $\sigma^\text{low}_i=h_{i,1} < h_{i,2}<...<h_{i,7}=\sigma^\text{high}_i$ and those values will be the ones populating the parameters $q_{e,j}$ according to the XBs' preferred days of the week for days off. Indeed, from their submitted days-off pattern preferences, we can compute the XBs' preferred days off from 1 to 7. In general, for each XB $e\in E$ and day $j\in J$ such that $j$ corresponds to the day of the week $i$, $q_{e,j}$ will be given the value $h_{i,n}$ where $e$ has day of the week $i$ as its $n^\text{th}$ least preferred day off. We illustrate this with an example on a division.
    
\begin{figure}[h!]
     \centerline{\includegraphics[width=\textwidth]{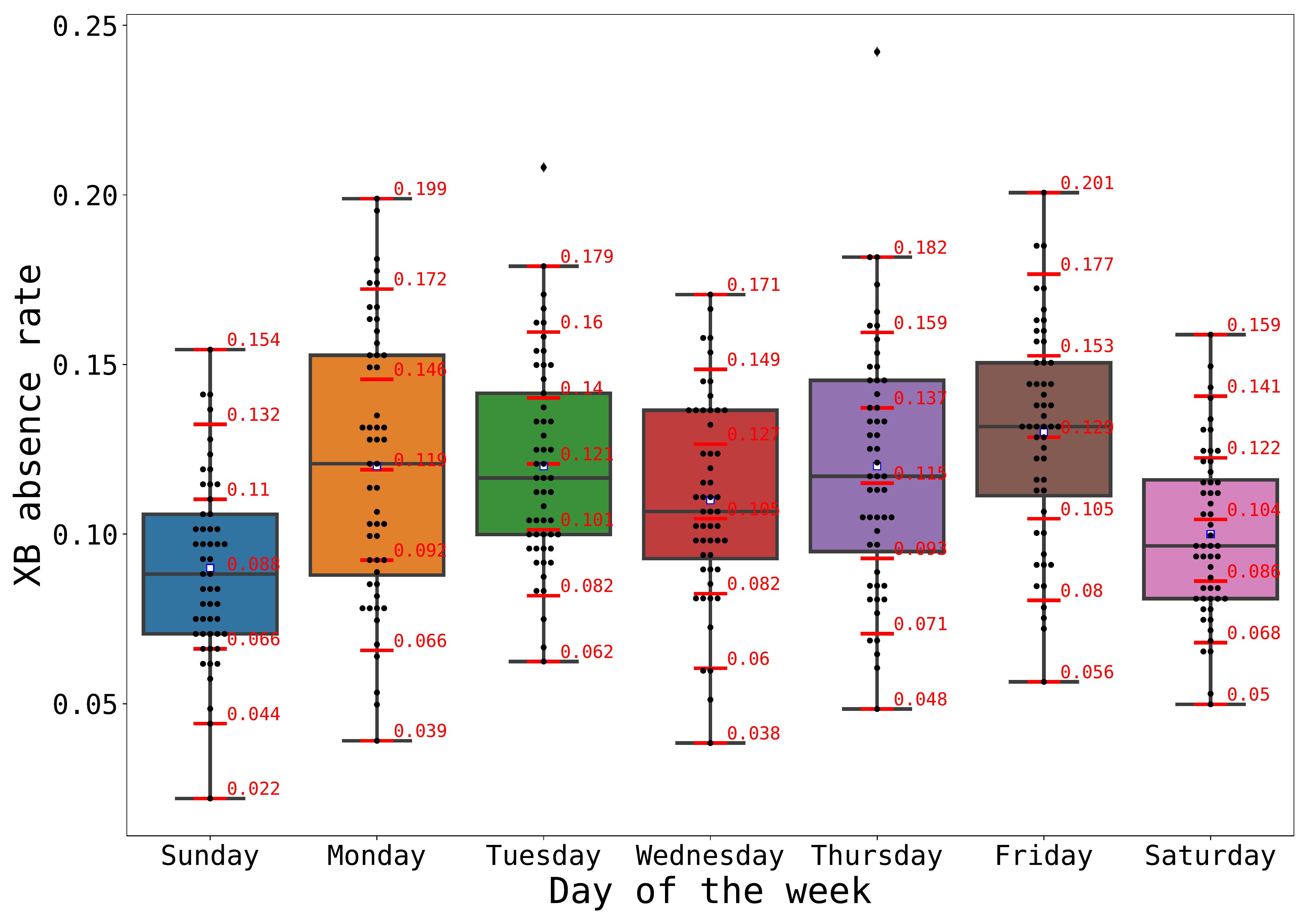}}
     \caption{Box plot of XB absence rates for division 1. For each day of the week $i$, the values $h_{i,1}$ to $h_{i,7}$ are displayed with the horizontal red lines.\label{fig:XBAbsProbs}}
     \end{figure}

    In \Cref{fig:XBAbsProbs}, we see that on Sundays ($i=1$), the low and high whiskers are $\sigma^\text{low}_1=0.022$ and $\sigma^\text{high}_1=0.154$, defining the range for $h_{1,1}, ..., h_{1,7}$, and on Wednesdays ($i=4$), $\sigma^\text{low}_4=0.038$ and $\sigma^\text{high}_4=0.171$. Then, suppose an XB $e$ submitted the following days-off pattern preferences (from most to least preferred): Sat-Sun (7), Sun-Mon (6), Fri-Sat (5), Thu-Fri (4), Wed-Thu (3), Tue-Wed (2), Mon-Tue (1). We are now able to rank his/her individual days of the week preferences for days off: Sun (13), Sat (12), Fri (9), Mon (7), Thu (7), Wed (5), Tue (3). In case of ties, the natural week ordering prevails. Then, on Wednesdays, XB $e$ has an absence probability of $h_{4,2}=0.06$ because Wednesdays are his/her $2^\text{nd}$ least desired day off. Hence, $q_{e,4}=q_{e,11}=0.06$ since days $j=4$ and $j=11$ of the horizon are Wednesdays. Similarly, on Sundays, his/her absence probability is of $h_{1,7}=0.154$. We can see in this example that the absence probability of XBs is high (low) in the whiskers interval, on days of the week where the XB has a strong (weak) desire to have a day off. It is also worth noting that by construction, for each day of the week $i\in \{1,2,...,7\}$, the values $q_i$ and $h_{i,4}$ are always very close to each other, making the absence probability of XBs receiving their middle choice very similar to the overall absence probability of XBs when their preferences are not considered.
     
\subsection{Other parameter settings}\label{sec:ParamSettings}
    In the experiments, we consider $|J|=14$-day long horizons, with an hourly discretization of the time ($|T|=24$).
    
    The number of XBs $|E|$ varies according to the division and can be found in \Cref{tab:StatsByDiv} along with an overview of the known-in-advance and unknown absence quantities. Due to a lack of information in the data, we proceeded to a sensitivity analysis to determine $|E|$ and fixed, for each division, the number of XBs minimizing the cost over the 25 planning horizons of the XB days-off problem without considering XB preferences.
    
    \begin{table}[h!]
    \caption{Number of XBs and overview of demand per division.\label{tab:StatsByDiv}}
    \begin{tabular*}{\hsize}{@{}@{\extracolsep{\fill}}ccccccccccc@{}}
    \hline
    \textbf{Division} & 1 & 2 & 3 & 4 & 5 & 6 & 7 & 8 & 9 & 10 \\ \hline
    $\bm{|E|}$ & 50 & 56 & 55 & 61 & 73 & 51 & 67 & 61 & 63 & 98 \\
    $\bm{\max_{j\in J}o_j}$ & 14 & 13 & 12 & 16 & 15 & 13 & 17 & 15 & 18 & 32 \\
    $\bm{\max_{j\in J} \text{LB}_j}$ & 30 & 31 & 35 & 36 & 45 & 32 & 39 & 34 & 40 & 46 \\ \hline
    \end{tabular*}
    \tablenote{{\it Note}: The third row indicates the average over all the planning horizons of the maximal daily number of known-in-advance absences over the horizon. $\text{LB}_j$ stands for a lower bound on the number of XBs to cover the unknown absences on day $j$, and is calculated by summing the daily unknown absence durations and dividing it by 8 hours (the duration of a regular XB shift). The last row indicates the average  over all scenarios and planning horizons of the maximal daily $\text{LB}_j$ over the horizon. The last two rows are rounded to nearest integer.}
    \end{table}
    
    The number of daily scenarios $|S_j|$ varies according to the experiments, but is always the same for all days of a horizon. When building these scenarios, the number of daily unknown absence durations $l$ is always equal to the number of unknown absence shapes $k$. The daily scenarios are always treated with equal probability: $\forall s_1, s_2\in S_j, \alpha^{s_1}_j=\alpha^{s_2}_j$.
    
    The set of duties $W$ contains 201 duties covering the whole day, with a potential 3 hours long pause in the middle (representing split duties, commonly used in public transit). Since the labor union agreement states that each employee must be paid a minimum of 8 hours, all the duties are created with a minimum of 8 work hours. After the first 8 work hours which are paid at regular rate, each additional working hour is paid as overtime. No duty exceeds a 12-hour span nor a 10-hour work duration, meaning that the maximum overtime duration is 2 hours. The cost $c^w_2$ of duties $w$ is calculated in the following way: a cost of $1.5$ is assessed for each overtime hour, $1$ for each regular working time hour and $0.5$ for each unpaid pause hour as a mental cost for the driver, making the cost of any duty ranging between 8 and 11.
    
    Costs $c_1=10$ and $c_3=0.75$ (when using XB preferences) were also tuned using a sensitivity analysis on a single division.
      
\subsection{Problem instances}\label{sec:ModelInstances}
    As stated in \Cref{sec:InputData}, we have in hand 25 planning horizons for each of the 10 divisions. We transform those planning horizons into 2500 problem instances by sampling 10 sets of XB days-off pattern preferences for each planning horizon and division. When not considering the days-off preferences, we remain with a total of 250 problem instances as we use the divisional average absence rates $q_j$ instead of the XB-preference related ones, as explained in \Cref{sec:ModelFormNoPref}.
    
    Whenever a problem instance is solved with XB preferences consideration, the first-stage solution is saved. It is then evaluated on a hidden set of 1000 second-stage problem instances by fixing the first-stage solution to the saved one. Each of these second-stage problems contains exactly 1 scenario from a held-out test set of 1000 scenarios.
    
    To ensure a fair comparison between solutions that do and do not take into account the XB days-off preferences, the solutions found in the latter case are evaluated on the 10 problem instances with different XB preference samplings, the test set of scenarios, and formulation (\ref{eq:ObjFun})-(\ref{eq:BinaryV}) instead of the special case formulation defined in \Cref{sec:ModelFormNoPref}. In this way, we can easily observe in the evaluation phase the repercussions of taking the days-off decisions without considering the XB preferences.

\subsection{Computational environment}\label{sec:CompEnv}
    The integer programs were all solved using the Gurobi MIP solver v9.5.0, limited to a single thread of an Intel Gold 6148 Skylake @ 2.4 GHz CPU. 16Gb of RAM were sufficient to solve all instances. In fact, the largest problem instances involve up to around 45,000 constraints and 365,000 integer variables. A time limit of 3 hours and an optimality gap tolerance of 0.01\% were given as stopping criteria.

\section{Computational experiments}\label{sec:CompExp}
    In this section, we report the experimental results. We start in \Cref{sec:ResultsPrefsOrNot} by showing the benefit of considering employee preferences. Then, in \Cref{sec:VSS}, we empirically motivate the interest of considering a stochastic formulation. \Cref{sec:NumScenarios,sec:EVPI} provide an analysis on the number of scenarios to use and a comparison with solutions under perfect information, respectively. The results shown in this section correspond to the evaluation phase described in \Cref{sec:ModelInstances}, and are averaged over the 1000 evaluation scenarios and the 25 planning horizons for each division.
    
\subsection{Impact of preferences}\label{sec:ResultsPrefsOrNot}
    \Cref{tab:Res_10-10-False-False,tab:Res_10-10-False-True} provide details about the stochastic solutions that do and do not consider the XB days-off preferences. \Cref{tab:Comp_table_with_without_prefs} compares these two sets of solutions.
    \newline
    
    \begin{table}[h!]
    \centering
    \caption{Statistics about stochastic solutions (100 daily scenarios) that do not consider preferences\label{tab:Res_10-10-False-False}}
    \begin{tabular*}{\hsize}{@{}@{\extracolsep{\fill}}ccccccccccc@{}}
    \hline
    \textbf{Division} & \textbf{Cost} & \textbf{\begin{tabular}[c]{@{}c@{}}C.S.\\ (h.)\end{tabular}} & \textbf{S.W.} & \textbf{XB abs.} & \textbf{\begin{tabular}[c]{@{}c@{}}OVS.\\ (h.)\end{tabular}} & \textbf{\begin{tabular}[c]{@{}c@{}}XB util.\\ rate\\ (\%)\end{tabular}} & \textbf{\begin{tabular}[c]{@{}c@{}}OVT.\\ (h.)\end{tabular}} & \textbf{\begin{tabular}[c]{@{}c@{}}Sol. time\\ (min.)\end{tabular}} & \textbf{\begin{tabular}[c]{@{}c@{}}Opt. gap\\ (\%)\end{tabular}} & \textbf{\begin{tabular}[c]{@{}c@{}}Solved\\ instances\end{tabular}} \\ \hline
    \textbf{1} & 2876.1 & 35.5 & 4.23 & \phantom{1}63.1 & \phantom{1}600.8 & 77.4 & 72.2 & 77.2 & 0.08 & 23/25 \\
    \textbf{2} & 3089.8 & 27.4 & 4.27 & \phantom{1}94.1 & \phantom{1}761.0 & 74.4 & 61.0 & 56.6 & 0.04 & 23/25 \\
    \textbf{3} & 3391.1 & 42.2 & 4.22 & \phantom{1}82.6 & \phantom{1}699.6 & 77.5 & 77.2 & 91.9 & 0.21 & 19/25 \\
    \textbf{4} & 3412.5 & 38.8 & 4.38 & 106.3           & \phantom{1}797.6 & 75.1 & 72.8 & 52.6 & 0.01 & 25/25 \\
    \textbf{5} & 4479.8 & 48.3 & 4.31 & 100.4           & \phantom{1}981.8 & 76.9 & 93.1 & 83.2 & 0.11 & 23/25 \\
    \textbf{6} & 3016.3 & 35.8 & 4.45 & \phantom{1}63.4 & \phantom{1}919.3 & 68.2 & 43.9 & 71.8 & 0.03 & 20/25 \\
    \textbf{7} & 3742.7 & 31.6 & 4.53 & \phantom{1}97.9 & \phantom{1}860.8 & 76.4 & 77.2 & 72.0 & 0.13 & 22/25 \\
    \textbf{8} & 3567.8 & 40.7 & 4.08 & \phantom{1}83.0 & \phantom{1}881.2 & 74.3 & 80.2 & 86.8 & 0.37 & 21/25 \\
    \textbf{9} & 3801.5 & 44.2 & 4.47 & \phantom{1}67.5 & 1013.2           & 72.1 & 55.0 & 49.4 & 0.18 & 23/25 \\
    \textbf{10} & 4761.6 & 41.7 & 4.48 & 153.0          & 1335.3           & 72.3 & 86.8 & 87.8 & 0.22 & 20/25 \\
    \textbf{Mean} & \textbf{3613.9} & \textbf{38.6} & \textbf{4.34} & \textbf{91.1} & \textbf{885.0} & \textbf{74.5} & \textbf{71.9} & \textbf{72.9} & \textbf{0.14} & \textbf{-} \\ \hline
    \end{tabular*}
    \tablenote{{\it Note}: The average and the maximal optimality gap for non-optimal instances are 1.06\% and 3.96\%, respectively.}
    \end{table}

    \begin{table}[h!]
    \centering
    \caption{Statistics about stochastic solutions (100 daily scenarios) that consider preferences\label{tab:Res_10-10-False-True}}
    \begin{tabular*}{\hsize}{@{}@{\extracolsep{\fill}}ccccccccccc@{}}
    \hline
    \textbf{Division} & \textbf{Cost} & \textbf{\begin{tabular}[c]{@{}c@{}}C.S.\\ (h.)\end{tabular}} & \textbf{S.W.} & \textbf{XB abs.} & \textbf{\begin{tabular}[c]{@{}c@{}}OVS.\\ (h.)\end{tabular}} & \textbf{\begin{tabular}[c]{@{}c@{}}XB util.\\ rate\\ (\%)\end{tabular}} & \textbf{\begin{tabular}[c]{@{}c@{}}OVT.\\ (h.)\end{tabular}} & \textbf{\begin{tabular}[c]{@{}c@{}}Sol. time\\ (min.)\end{tabular}} & \textbf{\begin{tabular}[c]{@{}c@{}}Opt. gap\\ (\%)\end{tabular}} & \textbf{\begin{tabular}[c]{@{}c@{}}Solved\\ instances\end{tabular}} \\ \hline
    \textbf{1} & 2773.0 & 28.7 & 5.86 & \phantom{1}58.5 & \phantom{1}623.3 & 76.7 & 64.7 & 180.0 & 2.12 & 0/250 \\
    \textbf{2} & 2988.8 & 21.6 & 5.78 & \phantom{1}89.7 & \phantom{1}781.2 & 73.9 & 51.5 & 180.0 & 1.69 & 0/250 \\
    \textbf{3} & 3276.1 & 34.1 & 5.85 & \phantom{1}76.5 & \phantom{1}730.6 & 76.7 & 67.5 & 180.0 & 1.75 & 0/250 \\
    \textbf{4} & 3270.4 & 28.5 & 5.89 & 100.7           & \phantom{1}821.7 & 74.5 & 62.7 & 180.0 & 1.57 & 0/250 \\
    \textbf{5} & 4327.7 & 39.3 & 6.00 & \phantom{1}94.9 & 1007.8           & 76.4 & 84.0 & 180.0 & 1.16 & 0/250 \\
    \textbf{6} & 2921.2 & 30.9 & 6.00 & \phantom{1}60.9 & \phantom{1}928.7 & 67.9 & 40.6 & 180.0 & 1.32 & 0/250 \\
    \textbf{7} & 3592.2 & 22.3 & 6.23 & \phantom{1}91.8 & \phantom{1}887.0 & 75.9 & 64.2 & 180.0 & 1.39 & 0/250 \\
    \textbf{8} & 3434.4 & 33.0 & 5.70 & \phantom{1}79.5 & \phantom{1}894.0 & 73.9 & 72.8 & 180.0 & 1.57 & 0/250 \\
    \textbf{9} & 3686.5 & 38.8 & 6.07 & \phantom{1}64.6 & 1025.9           & 71.8 & 50.1 & 180.0 & 1.19 & 0/250 \\
    \textbf{10} & 4547.9 & 28.2 & 6.04 & 145.6          & 1366.4           & 71.7 & 72.0 & 180.0 & 1.13 & 0/250 \\
    \textbf{Mean} & \textbf{3481.8} & \textbf{30.5} & \textbf{5.94} & \textbf{86.3} & \textbf{906.7} & \textbf{73.9} & \textbf{63.0} & \textbf{180.0} & \textbf{1.49} & \textbf{-} \\ \hline
    \end{tabular*}
    \tablenote{{\it Note}: The maximal optimality gap for non-optimal instances is 5.43\%.}
    \end{table}
    
    The cost, the cancelled service (C.S.), the average employee social welfare (S.W.), the number of XB absences (XB abs.), the overstaffing (OVS.), the XB utilization rate (XB util. rate, calculated by dividing the sum of hours where XBs cover unknown absences by the sum of their paid work hours) and the overtime (OVT.) are shown in \Cref{tab:Res_10-10-False-False,tab:Res_10-10-False-True}, and for similar tables in this section.  The solution time (Sol. time), the optimality gap (Opt. gap) and the number of instances solved to optimality (Solved instances) are also displayed and concern the first-stage solutions obtained from solving the extensive form of the two-stage problems (before the evaluation phase). Optimality is obtained in less than a second when solving any second-stage problem in the evaluation phase.
    \newline
    
    \begin{table}[h!]
    \centering
    \caption{Improvements from solutions that do not consider preferences to solutions that consider them}
    \label{tab:Comp_table_with_without_prefs}
    \begin{tabular}{cccc}
    \hline
    \textbf{Division} & \textbf{\begin{tabular}[c]{@{}c@{}}Cost\\ (\%)\end{tabular}} & \textbf{\begin{tabular}[c]{@{}c@{}}C.S.\\ (\%)\end{tabular}} & \textbf{\begin{tabular}[c]{@{}c@{}}S.W.\\ (\%)\end{tabular}} \\ \hline
    \textbf{1} & 3.34 & 19.01 & 38.87 \\
    \textbf{2} & 3.11 & 21.09 & 35.70 \\
    \textbf{3} & 3.13 & 19.35 & 39.22 \\
    \textbf{4} & 3.82 & 26.49 & 34.64 \\
    \textbf{5} & 3.21 & 18.72 & 39.43 \\
    \textbf{6} & 3.04 & 13.58 & 35.43 \\
    \textbf{7} & 3.88 & 29.60 & 37.72 \\
    \textbf{8} & 3.54 & 18.81 & 40.03 \\
    \textbf{9} & 2.76 & 12.39 & 36.18 \\
    \textbf{10} & 4.39 & 32.34 & 34.92 \\
    \textbf{Mean} & \textbf{3.42} & \textbf{21.14} & \textbf{37.21} \\ \hline
    \end{tabular}
    \end{table}
    
    In \Cref{tab:Comp_table_with_without_prefs} and for similar tables in the rest of the section, the average over all evaluation instances of the cost and social welfare ratios over the two methods are shown. For the cancelled service, the ratio of the average are displayed to account for solutions with no cancelled service on some evaluation instances.
    
    We observe important differences between the two approaches. The solutions with preferences outclass the ones without preferences in terms of cancelled service, social welfare and cost. In fact, the solutions with preferences are on average close to a social welfare of 6, meaning that the XBs are assigned on average to their $2^\text{nd}$ most preferred days-off pattern. Without preferences, the social welfare reaches an average value of 4.34 (close to the $4^\text{th}$ most preferred pattern). This 37.21\% social welfare improvement is directly reflected on the level of employee satisfaction which permits to save on average around 5 XB absences per planning horizon when considering preferences. Consequently, the cancelled service and the total cost are much lower when considering preferences, beating the solutions without preferences by 21.14\% and 3.42\% respectively. Other consequences of the 5 XB absences saved are the increase in overstaffing and decrease in utilization rate from the solutions that do not consider preferences to the ones that do; indeed, the additional XB duties that we can assign might not cover a lot of demand, but help to cover the demand peaks during the horizon. The overtime also follows a decrease, which is intuitive. In terms of solution time, however, the stochastic program with preferences is never solved to proven optimality within 3 hours, while the programs without preferences are solved on average in 1h13min.

\subsection{Value of the stochastic solution}\label{sec:VSS}
    The value of stochastic solution (VSS) is introduced in \cite{Birge1982}, and expresses how well the stochastic formulation takes care of the problem uncertainty compared to the deterministic formulation. The deterministic formulation considers only 1 daily scenario each day, that we choose to be the average of the 100 scenarios from the stochastic formulation. Details of the solutions obtained by ignoring and considering preferences are shown in \Cref{tab:10-10-True-False,tab:10-10-True-True}, respectively. The VSSs are reported in \Cref{tab:Comp_table_VSS}.
    \newline
    
    \begin{table}[h!]
    \centering
    \caption{Statistics about deterministic solutions that do not consider preferences}
    \label{tab:10-10-True-False}
    \begin{tabular*}{\hsize}{@{}@{\extracolsep{\fill}}ccccccccccc@{}}
    \hline
    \textbf{Division} & \textbf{Cost} & \textbf{\begin{tabular}[c]{@{}c@{}}C.S.\\ (h.)\end{tabular}} & \textbf{S.W.} & \textbf{XB abs.} & \textbf{\begin{tabular}[c]{@{}c@{}}OVS.\\ (h.)\end{tabular}} & \textbf{\begin{tabular}[c]{@{}c@{}}XB util.\\ rate\\ (\%)\end{tabular}} & \textbf{\begin{tabular}[c]{@{}c@{}}OVT.\\ (h.)\end{tabular}} & \textbf{\begin{tabular}[c]{@{}c@{}}Sol. time\\ (min.)\end{tabular}} & \textbf{\begin{tabular}[c]{@{}c@{}}Opt. gap\\ (\%)\end{tabular}} & \textbf{\begin{tabular}[c]{@{}c@{}}Solved\\ instances\end{tabular}} \\ \hline
    \textbf{1} & 2932.3 & 40.3 & 4.29 & \phantom{1}62.6 & \phantom{1}613.9 & 77.0 & 76.7 & 0.0 & 0.00 & 25/25 \\
    \textbf{2} & 3170.2 & 35.5 & 4.29 & \phantom{1}94.5 & \phantom{1}771.2 & 74.1 & 66.0 & 0.0 & 0.00 & 25/25 \\
    \textbf{3} & 3435.6 & 46.3 & 4.33 & \phantom{1}82.0 & \phantom{1}711.7 & 77.2 & 81.0 & 0.0 & 0.00 & 25/25 \\
    \textbf{4} & 3452.0 & 42.2 & 4.44 & 106.1           & \phantom{1}807.9 & 74.8 & 78.5 & 0.0 & 0.00 & 25/25 \\
    \textbf{5} & 4501.7 & 50.0 & 4.37 & 100.3           & \phantom{1}990.2 & 76.8 & 98.9 & 0.0 & 0.00 & 25/25 \\
    \textbf{6} & 3050.7 & 38.9 & 4.59 & \phantom{1}63.1 & \phantom{1}929.7 & 68.0 & 48.5 & 0.1 & 0.00 & 25/25 \\
    \textbf{7} & 3786.8 & 35.4 & 4.58 & \phantom{1}97.6 & \phantom{1}871.3 & 76.2 & 81.6 & 0.0 & 0.00 & 25/25 \\
    \textbf{8} & 3592.7 & 42.9 & 4.00 & \phantom{1}83.6 & \phantom{1}881.3 & 74.3 & 83.5 & 0.6 & 0.00 & 25/25 \\
    \textbf{9} & 3824.4 & 46.1 & 4.50 & \phantom{1}66.9 & 1020.8           & 71.9 & 56.7 & 0.2 & 0.00 & 25/25 \\
    \textbf{10} & 4796.7 & 44.0 & 4.54 & 152.6          & 1348.7           & 72.1 & 94.5 & 0.2 & 0.00 & 25/25 \\
    \textbf{Mean} & \textbf{3654.3} & \textbf{42.2} & \textbf{4.39} & \textbf{\phantom{1}90.9} & \textbf{\phantom{1}894.7} & \textbf{74.2} & \textbf{76.6} & \textbf{0.1} & \textbf{0.00} & \textbf{-} \\ \hline
    \end{tabular*}
    \end{table}
    
    \begin{table}[h!]
    \centering
    \caption{Statistics about deterministic solutions that consider preferences}
    \label{tab:10-10-True-True}
    \begin{tabular*}{\hsize}{@{}@{\extracolsep{\fill}}ccccccccccc@{}}
    \hline
    \textbf{Division} & \textbf{Cost} & \textbf{\begin{tabular}[c]{@{}c@{}}C.S.\\ (h.)\end{tabular}} & \textbf{S.W.} & \textbf{XB abs.} & \textbf{\begin{tabular}[c]{@{}c@{}}OVS.\\ (h.)\end{tabular}} & \textbf{\begin{tabular}[c]{@{}c@{}}XB util.\\ rate\\ (\%)\end{tabular}} & \textbf{\begin{tabular}[c]{@{}c@{}}OVT.\\ (h.)\end{tabular}} & \textbf{\begin{tabular}[c]{@{}c@{}}Sol. time\\ (min.)\end{tabular}} & \textbf{\begin{tabular}[c]{@{}c@{}}Opt. gap\\ (\%)\end{tabular}} & \textbf{\begin{tabular}[c]{@{}c@{}}Solved\\ instances\end{tabular}} \\ \hline
    \textbf{1} & 2810.4 & 35.2 & 5.68 & \phantom{1}64.5 & \phantom{1}590.9 & 77.6 & 73.7 & 161.2 & 0.28 & 36/250 \\
    \textbf{2} & 3042.8 & 29.6 & 5.16 & \phantom{1}98.5 & \phantom{1}731.0 & 75.1 & 63.9 & 159.6 & 0.20 & 41/250 \\
    \textbf{3} & 3314.1 & 41.6 & 5.25 & \phantom{1}86.1 & \phantom{1}669.9 & 78.2 & 76.6 & 172.0 & 0.16 & 16/250 \\
    \textbf{4} & 3312.7 & 35.4 & 5.05 & 111.4           & \phantom{1}755.3 & 76.0 & 74.8 & 168.0 & 0.23 & 30/250 \\
    \textbf{5} & 4353.8 & 43.5 & 5.72 & 101.2           & \phantom{1}972.5 & 77.0 & 95.2 & 174.5 & 0.22 & 12/250 \\
    \textbf{6} & 2938.1 & 34.4 & 5.69 & \phantom{1}65.6 & \phantom{1}897.8 & 68.7 & 45.1 & 146.7 & 0.14 & 68/250 \\
    \textbf{7} & 3661.6 & 31.2 & 5.75 & 100.2           & \phantom{1}842.5 & 76.8 & 78.3 & 170.6 & 0.18 & 21/250 \\
    \textbf{8} & 3456.4 & 37.6 & 5.36 & \phantom{1}86.1 & \phantom{1}854.2 & 74.8 & 81.0 & 179.7 & 0.32 & 1/250 \\
    \textbf{9} & 3696.5 & 41.2 & 5.99 & \phantom{1}68.5 & 1004.2           & 72.2 & 57.2 & 166.1 & 0.23 & 29/250 \\
    \textbf{10} & 4593.0 & 33.8 & 5.87 & 152.1          & 1336.6           & 72.3 & 88.1 & 180.0 & 0.27 & 0/250 \\
    \textbf{Mean} & \textbf{3517.9} & \textbf{36.4} & \textbf{5.55} & \textbf{\phantom{1}93.4} & \textbf{\phantom{1}865.5} & \textbf{74.9} & \textbf{73.4} & \textbf{167.8} & \textbf{0.22} & \textbf{-} \\ \hline
    \end{tabular*}
    \tablenote{{\it Note}: The average and the maximal optimality gap for non-optimal instances are 0.25\% and 2.87\%, respectively.}
    \end{table}
    
    \begin{table}[h!]
    \centering
    \caption{Value of stochastic solutions: improvements from deterministic solutions to stochastic solutions using 100 daily scenarios\label{tab:Comp_table_VSS}}
    \begin{tabular}{cccccccc}
    \hline
     & \multicolumn{3}{c}{\textbf{\begin{tabular}[c]{@{}c@{}}VSS\\ without preferences\end{tabular}}} & \textbf{} & \multicolumn{3}{c}{\textbf{\begin{tabular}[c]{@{}c@{}}VSS\\ with preferences\end{tabular}}} \\ \cline{2-4} \cline{6-8} 
    \textbf{Division} & \textbf{\begin{tabular}[c]{@{}c@{}}Cost\\ (\%)\end{tabular}} & \textbf{\begin{tabular}[c]{@{}c@{}}C.S.\\ (\%)\end{tabular}} & \textbf{\begin{tabular}[c]{@{}c@{}}S.W.\\ (\%)\end{tabular}} &  & \textbf{\begin{tabular}[c]{@{}c@{}}Cost\\ (\%)\end{tabular}} & \textbf{\begin{tabular}[c]{@{}c@{}}C.S.\\ (\%)\end{tabular}} & \textbf{\begin{tabular}[c]{@{}c@{}}S.W.\\ (\%)\end{tabular}} \\ \hline
    \textbf{1} & 1.71 & 12.10 & -1.02 &  & 1.10 & 18.33 & \phantom{1}3.31 \\
    \textbf{2} & 2.28 & 22.72 & -0.16 &  & 1.58 & 26.96 & 12.62 \\
    \textbf{3} & 1.11 & \phantom{1}8.78 & -2.22 &  & 0.87 & 18.15 & 12.18 \\
    \textbf{4} & 1.06 & \phantom{1}8.12 & -1.14 &  & 1.14 & 19.39 & 17.88 \\
    \textbf{5} & 0.41 & \phantom{1}3.34 & -1.15 &  & 0.47 & \phantom{1}9.80 & \phantom{1}5.02 \\
    \textbf{6} & 1.02 & \phantom{1}8.11 & -2.90 &  & 0.42 & 10.28 & \phantom{1}5.71 \\
    \textbf{7} & 1.09 & 10.58 & -0.75 &  & 1.73 & 28.61 & \phantom{1}8.64 \\
    \textbf{8} & 0.54 & \phantom{1}5.04 & \phantom{-}2.51 &  & 0.48 & 12.12 & \phantom{1}6.62 \\
    \textbf{9} & 0.61 & \phantom{1}3.95 & -0.27 &  & 0.25 & \phantom{1}6.00 & \phantom{1}1.56 \\
    \textbf{10} & 0.62 & \phantom{1}5.30 & -1.28 &  & 0.88 & 16.65 & \phantom{1}2.97 \\
    \textbf{Mean} & \textbf{1.04} & \textbf{\phantom{1}8.80} & \textbf{-0.84} &  & \textbf{0.89} & \textbf{16.63} & \textbf{\phantom{1}7.65} \\ \hline
    \end{tabular}
    \end{table}
    
    The behaviours are different when comparing the stochastic solutions to the deterministic solutions, depending on whether the preferences are taken into account or not. Although the VSSs are similar on average (1.04\% by not considering preferences, 0.89\% by considering preferences), it is when looking at the social welfare and cancelled service improvements that we see evidence of the strength of the stochastic formulation. The cancelled service is reduced by 8.80\% and 16.63\% (by not considering and by considering the preferences, respectively) when using the stochastic solutions over deterministic ones. Similarly, when considering the preferences, the social welfare is subject to an increase of 7.65\%. Given that social welfare maximization is not considered in the stochastic and deterministic solutions without preferences, it is not surprising that no social welfare trend is observed when comparing the respective solutions. We encounter the same phenomenon explained earlier about the reduction of XB absences, leading to a decrease in utilization rate and overtime, and an increase in overstaffing from the stochastic solutions to the deterministic solutions, when preferences are taken into account. In terms of solution time, it is interesting to notice that most of the first-stage solutions for the deterministic formulation with preferences are not proven to be optimal. Recall that the deterministic formulation considers only 1 scenario, highlighting the difficulty of the problem.
    
    Counter-intuitively, we can see that for the deterministic formulations, although the social welfare is greater when preferences are considered compared to when they are not considered, the number of XB absences is also superior. By looking at the details of the solutions, we explain this by the fact that the deterministic model with preferences saves the cost of a few duties by allowing more absences to occur on days where it is easy to cover the demand. Increasing the number of XB absences cannot be done without decreasing the social welfare. The maximal social welfare decrease for an XB $e$ costs at most $(\max_{p\in P}s_{p,e}-\min_{p\in P}s_{p,e})\cdot c_3=(7-1)\cdot 0.75=4.5$, and the cost of saving a duty is at least $\min_{w\in W}c_2^w=8$. Hence, the model tends to make few XBs absent when it does not result in cancelled service. We note that this behaviour disappears once multiple scenarios are taken into consideration. This is because, in that setting, having one more XB absence could create some cancelled service on some specific scenarios, inferring a cost of $c_1=10$ for these scenarios and making the duty cost savings not worth anymore.
    
\subsection{Varying the number of scenarios}\label{sec:NumScenarios}
    In this section, we justify our choice of using 100 daily scenarios. \Cref{tab:Comp_table_scenarios_49_25} compares the solutions obtained when using 25 and 49 daily scenarios, and \Cref{tab:Comp_table_scenarios_100_49}, for 49 and 100 scenarios.
    \newline

    \begin{table}[h!]
    \centering
    \caption{Improvements from solutions using 25 daily scenarios to solutions using 49 daily scenarios}
    \label{tab:Comp_table_scenarios_49_25}
    \begin{tabular}{cccccccc}
    \hline
     & \multicolumn{3}{c}{\textbf{\begin{tabular}[c]{@{}c@{}}49 Vs 25 scenarios\\ without preferences\end{tabular}}} & \textbf{} & \multicolumn{3}{c}{\textbf{\begin{tabular}[c]{@{}c@{}}49 Vs 25 scenarios\\ with preferences\end{tabular}}} \\ \cline{2-4} \cline{6-8} 
    \textbf{Division} & \textbf{\begin{tabular}[c]{@{}c@{}}Cost\\ (\%)\end{tabular}} & \textbf{\begin{tabular}[c]{@{}c@{}}C.S.\\ (\%)\end{tabular}} & \textbf{\begin{tabular}[c]{@{}c@{}}S.W.\\ (\%)\end{tabular}} &  & \textbf{\begin{tabular}[c]{@{}c@{}}Cost\\ (\%)\end{tabular}} & \textbf{\begin{tabular}[c]{@{}c@{}}C.S.\\ (\%)\end{tabular}} & \textbf{\begin{tabular}[c]{@{}c@{}}S.W.\\ (\%)\end{tabular}} \\ \hline
    \textbf{1} & \phantom{-}0.98 & \phantom{-}7.93 & \phantom{-}0.44 &  & 0.84 & 8.25 & \phantom{-}0.73 \\
    \textbf{2} & \phantom{-}0.75 & \phantom{-}7.29 & \phantom{-}0.83 &  & 0.49 & 6.20 & \phantom{-}0.39 \\
    \textbf{3} & \phantom{-}0.82 & \phantom{-}6.68 & \phantom{-}0.63 &  & 0.55 & 5.56 & \phantom{-}0.54 \\
    \textbf{4} & \phantom{-}0.55 & \phantom{-}5.20 & -0.27 &  & 0.49 & 4.54 & \phantom{-}0.22 \\
    \textbf{5} & \phantom{-}0.47 & \phantom{-}3.92 & \phantom{-}0.12 &  & 0.32 & 2.83 & \phantom{-}0.44 \\
    \textbf{6} & \phantom{-}0.20 & \phantom{-}1.73 & -1.51 &  & 0.30 & 3.20 & -0.14 \\
    \textbf{7} & -0.20 & -1.31 & \phantom{-}0.02 &  & 0.19 & 3.24 & \phantom{-}0.36 \\
    \textbf{8} & \phantom{-}0.67 & \phantom{-}6.12 & \phantom{-}0.73 &  & 0.65 & 7.60 & \phantom{-}0.72 \\
    \textbf{9} & \phantom{-}0.17 & \phantom{-}1.42 & \phantom{-}0.29 &  & 0.06 & 0.63 & \phantom{-}0.05 \\
    \textbf{10} & \phantom{-}0.70 & \phantom{-}6.95 & -0.41 &  & 0.22 & 3.85 & \phantom{-}0.26 \\
    \textbf{Mean} & \textbf{\phantom{-}0.51} & \textbf{\phantom{-}4.59} & \textbf{\phantom{-}0.09} &  & \textbf{0.41} & \textbf{4.59} & \textbf{\phantom{-}0.36} \\ \hline
    \tablenote{{\it Note}: When using 25 daily scenarios and without preferences, the average solution time is 17.3min and 241/250 instances are solved optimally. With preferences, none of the 2500 are solved optimally in 3h; the average and maximal optimality gaps are 1.19\% and 4.45\%, respectively.}
    \end{tabular}
    \end{table}
    
    \begin{table}[h!]
    \centering
    \caption{Improvements from solutions using 49 daily scenarios to solutions using 100 daily scenarios}
    \label{tab:Comp_table_scenarios_100_49}
    \begin{tabular}{cccccccc}
    \hline
     & \multicolumn{3}{c}{\textbf{\begin{tabular}[c]{@{}c@{}}100 Vs 49 scenarios\\ without preferences\end{tabular}}} & \textbf{} & \multicolumn{3}{c}{\textbf{\begin{tabular}[c]{@{}c@{}}100 Vs 49 scenarios\\ with preferences\end{tabular}}} \\ \cline{2-4} \cline{6-8} 
    \textbf{Division} & \textbf{\begin{tabular}[c]{@{}c@{}}Cost\\ (\%)\end{tabular}} & \textbf{\begin{tabular}[c]{@{}c@{}}C.S.\\ (\%)\end{tabular}} & \textbf{\begin{tabular}[c]{@{}c@{}}S.W.\\ (\%)\end{tabular}} &  & \textbf{\begin{tabular}[c]{@{}c@{}}Cost\\ (\%)\end{tabular}} & \textbf{\begin{tabular}[c]{@{}c@{}}C.S.\\ (\%)\end{tabular}} & \textbf{\begin{tabular}[c]{@{}c@{}}S.W.\\ (\%)\end{tabular}} \\ \hline
    \textbf{1} & 0.38 & \phantom{-}3.39 & \phantom{-}0.94 &  & 0.34 & 4.01 & \phantom{-}0.41 \\
    \textbf{2} & 0.16 & \phantom{-}1.40 & \phantom{-}0.08 &  & 0.30 & 4.68 & \phantom{-}0.75 \\
    \textbf{3} & 0.38 & \phantom{-}2.37 & \phantom{-}0.06 &  & 0.43 & 3.57 & \phantom{-}0.93 \\
    \textbf{4} & 0.27 & \phantom{-}1.16 & \phantom{-}0.55 &  & 0.33 & 4.32 & \phantom{-}0.98 \\
    \textbf{5} & 0.08 & -0.01 & \phantom{-}0.40 &  & 0.02 & 1.18 & -0.18 \\
    \textbf{6} & 0.26 & \phantom{-}2.02 & \phantom{-}0.95 &  & 0.10 & 1.66 & \phantom{-}0.71 \\
    \textbf{7} & 0.36 & \phantom{-}3.71 & -0.26 &  & 0.18 & 3.96 & \phantom{-}0.45 \\
    \textbf{8} & 0.38 & \phantom{-}3.45 & \phantom{-}0.54 &  & 0.31 & 3.23 & \phantom{-}0.38 \\
    \textbf{9} & 0.27 & \phantom{-}3.49 & \phantom{-}0.77 &  & 0.08 & 1.11 & \phantom{-}0.68 \\
    \textbf{10} & 0.10 & \phantom{-}1.13 & \phantom{-}1.01 &  & 0.16 & 1.96 & \phantom{-}0.93 \\
    \textbf{Mean} & \textbf{0.26} & \textbf{\phantom{-}2.21} & \textbf{\phantom{-}0.50} &  & \textbf{0.22} & \textbf{2.97} & \textbf{\phantom{-}0.60} \\ \hline
    \tablenote{{\it Note}: When using 49 daily scenarios and without preferences, the average solution time is 36.1min and 242/250 instances are solved optimally. With preferences, none of the 2500 are solved optimally in 3h; the average and maximal optimality gaps are 1.34\% and 5.42\%, respectively.}
    \end{tabular}
    \end{table}
    
    When switching from 25 to 49 daily scenarios, the total cost is improved on average by 0.51\% and 0.41\%, for the case without and with preferences, respectively. These values drop to 0.26\% and 0.22\% when going from 49 to 100 daily scenarios. When taking into account the preferences, the average optimality gap is slowly increasing, going from 1.19\% when using 25 daily scenarios to 1.34\% with 49 scenarios and finally 1.49\% with 100 scenarios. The slow increase in total cost improvements along with the increase in optimality gaps confirms the sufficiency of using 100 daily scenarios.
    
    It is also worth mentioning that the cancelled service is much more sensitive to the number of scenarios than the social welfare. This is explained by the direct relationship between the cancelled service and the demand that the scenarios represent. The relation between scenarios and social welfare is less obvious, making its improvements small when increasing the number of scenarios.

\subsection{Expected value of perfect information}\label{sec:EVPI}
    The expected value of prefect information (EVPI) is introduced in \cite{Pratt95} and measures the loss due to the problem uncertainty. It is calculated by comparing the stochastic and the perfect information solution costs on each evaluation scenario, where the perfect information solution costs are obtained by solving the problem after having observed the realization of the random variables (also called wait-and-see approach). Due to the high computational time to find the perfect information solutions that consider preferences on all the problem instances (3h for each of the 2,500,000 instances), we restrict ourselves to the case with no preferences. \Cref{tab:EVPI-False} shows the details from the perfect information solutions without considering preferences. The EVPI are then reported in \Cref{tab:Comp_table_EVPI}. 
    \newline

    \begin{table}[h!]
    \centering
    \caption{Statistics about perfect information solutions that do not consider preferences}
    \label{tab:EVPI-False}
    \begin{tabular*}{\hsize}{@{}@{\extracolsep{\fill}}ccccccccccc@{}}
    \hline
    \textbf{Division} & \textbf{Cost} & \textbf{\begin{tabular}[c]{@{}c@{}}C.S.\\ (h.)\end{tabular}} & \textbf{S.W.} & \textbf{XB abs.} & \textbf{\begin{tabular}[c]{@{}c@{}}OVS.\\ (h.)\end{tabular}} & \textbf{\begin{tabular}[c]{@{}c@{}}XB util.\\ rate\\ (\%)\end{tabular}} & \textbf{\begin{tabular}[c]{@{}c@{}}OVT.\\ (h.)\end{tabular}} & \textbf{\begin{tabular}[c]{@{}c@{}}Sol. time\\ (min.)\end{tabular}} & \textbf{\begin{tabular}[c]{@{}c@{}}Opt. gap\\ (\%)\end{tabular}} & \textbf{\begin{tabular}[c]{@{}c@{}}Solved\\ instances\end{tabular}} \\ \hline
    \textbf{1} & 2654.8 & 15.0 & 4.25 & \phantom{1}62.7 & \phantom{1}569.7 & 78.5 & 58.3 & 0.0 & 0.00 & 249820/250000 \\
    \textbf{2} & 2927.0 & 13.1 & 4.27 & \phantom{1}94.0 & \phantom{1}736.7 & 75.2 & 49.6 & 0.0 & 0.00 & 250000/250000 \\
    \textbf{3} & 3127.5 & 17.6 & 4.30 & \phantom{1}82.0 & \phantom{1}667.3 & 78.5 & 64.8 & 0.0 & 0.00 & 250000/250000 \\
    \textbf{4} & 3203.4 & 19.9 & 4.43 &           106.0 & \phantom{1}768.1 & 75.9 & 60.4 & 0.0 & 0.00 & 249980/250000 \\
    \textbf{5} & 4188.0 & 22.3 & 4.31 &           100.4 & \phantom{1}933.8 & 78.0 & 71.3 & 0.0 & 0.00 & 249880/250000 \\
    \textbf{6} & 2876.3 & 23.5 & 4.50 & \phantom{1}63.4 & \phantom{1}898.3 & 68.9 & 34.6 & 0.0 & 0.00 & 244290/250000 \\
    \textbf{7} & 3517.0 & 12.1 & 4.54 & \phantom{1}97.7 & \phantom{1}822.8 & 77.4 & 57.4 & 0.0 & 0.00 & 249920/250000 \\
    \textbf{8} & 3306.2 & 17.5 & 4.05 & \phantom{1}83.4 & \phantom{1}836.3 & 75.5 & 62.1 & 0.0 & 0.00 & 249980/250000 \\
    \textbf{9} & 3673.3 & 34.2 & 4.42 & \phantom{1}67.5 & \phantom{1}985.3 & 72.8 & 37.3 & 0.0 & 0.00 & 249860/250000 \\
    \textbf{10} & 4513.4 & 19.1 & 4.48 &          153.1 &           1298.4 & 73.1 & 72.7 & 0.0 & 0.00 & 249640/250000 \\
    \textbf{Mean} & \textbf{3398.7} & \textbf{19.4} & \textbf{4.35} & \textbf{\phantom{1}91.0} & \textbf{\phantom{1}851.7} & \textbf{75.4} & \textbf{56.9} & \textbf{0.0} & \textbf{0.00} & \textbf{-} \\ \hline
    \end{tabular*}
    \tablenote{{\it Note}: The average and the maximal optimality gap for non-optimal instances are 0.11\% and 0.44\%, respectively.}
    \end{table}
    
    \begin{table}[h!]
    \centering
    \caption{Expected value of perfect information that do not consider preferences: improvements from stochastic solutions using 100 daily scenarios to perfect information solutions}
    \label{tab:Comp_table_EVPI}
    \begin{tabular}{cccc}
    \hline
    \textbf{Division} & \textbf{\begin{tabular}[c]{@{}c@{}}Cost\\ (\%)\end{tabular}} & \textbf{\begin{tabular}[c]{@{}c@{}}C.S.\\ (\%)\end{tabular}} & \textbf{\begin{tabular}[c]{@{}c@{}}S.W.\\ (\%)\end{tabular}} \\ \hline
    \textbf{1} & 7.20 & 57.63 & \phantom{-}0.67 \\
    \textbf{2} & 5.01 & 52.05 & \phantom{-}0.41 \\
    \textbf{3} & 7.03 & 58.28 & \phantom{-}2.28 \\
    \textbf{4} & 5.72 & 48.61 & \phantom{-}1.36 \\
    \textbf{5} & 6.04 & 53.88 & \phantom{-}0.10 \\
    \textbf{6} & 4.41 & 34.18 & \phantom{-}1.46 \\
    \textbf{7} & 5.76 & 61.72 & \phantom{-}0.38 \\
    \textbf{8} & 6.85 & 56.99 & -0.44 \\
    \textbf{9} & 3.20 & 22.59 & -0.94 \\
    \textbf{10} & 5.08 & 54.30 & \phantom{-}0.29 \\
    \textbf{Mean} & \textbf{5.63} & \textbf{50.02} & \textbf{\phantom{-}0.56} \\ \hline
    \end{tabular}
    \end{table}
    
    Since the perfect information solutions without preferences are obtained by solving the formulation for the case without preferences in \Cref{sec:ModelFormNoPref} and then evaluated on the second-stage evaluation problem instances with XB preferences, we can sometimes witness the stochastic solutions without preferences outperforming the perfect information solutions. As a matter of fact, 227,111 stochastic solutions obtained a lower cost than the perfect information solutions out of the 2,500,000 evaluation instances, representing around 9.08\%. Also, because solving the deterministic formulation without preferences, it is not surprising to see the average solution time and optimality gap of 0.0min and 0.00\%, respectively, for the perfect information solutions.
    
    The EVPI is already high at 5.63\% on average, but we can particularly see the effect of the uncertainty on the cancelled service: knowing in advance the demand would allow us to decrease on average by one half the cancelled service. We expect these values to be even higher when considering preferences. 

\section{Conclusion}\label{sec:Conc}

    In this paper, the extra board days-off scheduling problem is studied for the first time. We develop a flexible two-stage stochastic integer program that can take into consideration the employee preferences. This problem is interesting because of its applicability in the transit industry. From the practical side, an improvement would be to include an additional recourse action, which is to call regular drivers to work overtime at the end of their duties or during their days off. This action is regularly used in the transit companies because it permits to save service punctually without hiring additional extra boards that would decrease the global extra board utilization rate.
    
    The problem is also methodologically challenging. Indeed, the addition of endogenous uncertainty related to the employee preferences and absence probabilities, to the original exogenous demand uncertainty, leaves unsolved even the smallest mixed-integer program instances that involve as low as 500 constraints and 4000 variables. Nevertheless, based on instances generated from real-world data, we show empirically the added value of representing uncertainty through a stochastic program as well as of considering employee satisfaction to decrease their absences and the cancelled service. From the methodological side, a great addition to our work would be to derive a solution approach to solve our mixed-integer programs to optimality as further improvements on cost, cancelled service and social welfare may be obtained.
    
    Our work is currently being implemented in a world-leader software solutions provider for public transit operations.

\section*{Ethical Impact}\label{sec:EthicalImpact}
    The pipeline we develop in this work shows that taking extra board preferences into account has benefits for both the transit company and the employees. The cancelled service reduction is greater and the extra board social welfare is higher than when preferences are not considered. However, particular attention should be given to the individual preference scores of the days-off that extra boards receive. Indeed, since we are globally maximizing the social welfare, nothing prevents the dispersion of values forming the welfare and some extra boards could receive days off they do not like even if the social welfare is high. Any real-world implementation of this work should carefully monitor the individual preference scores to avoid unlucky XBs receiving disliked days off, horizon after horizon.

\section*{Acknowledgments}\label{sec:Acknowledgments}
    We would like to thank the GIRO Inc. team for their support and the Los Angeles County Metropolitan Transportation Authority for sharing their data. This work was funded by the FRQ-IVADO Research Chair in Data Science for Combinatorial Game Theory and the NSERC discovery grant 2019-04557. We also wish to thank Hydro-Québec for their financial support via the 2021 Excellence Scholarships. This research was enabled in part by support provided by Calcul Québec (\url{www.calculquebec.ca}) and Digital Research Alliance of Canada (\url{www.alliancecan.ca})

\nocite{*}

\bibliographystyle{itor}
\bibliography{itor}

\end{document}